\documentclass[]{article}
\usepackage{geometry}
\usepackage{amsmath,amssymb,amsfonts,amsthm,amsopn,dsfont, mathrsfs}
\usepackage{bm}
\usepackage{graphicx}
\usepackage{float}

\usepackage{standalone}
\usepackage{tikz}
\usepackage{pgfplots}
\usepackage{pgfplotstable}
\pgfplotsset{compat=newest}
\usepackage{subcaption}
\usepackage{colortbl}
\usepackage{booktabs}
\usepackage[backend=bibtex,citestyle=numeric-comp]{biblatex}

\usepackage[T1]{fontenc}

\usepackage{xcolor}

\usepackage{multirow}

\newcommand{\flux}{\bm{\sigma}_h}
\newcommand{\fluxf}{\tilde{\bm{\sigma}}_h}
\newcommand{\fluxa}{\bm{\sigma}_h^{\bm{a}}}

\newcommand{\tr}[2]{{#1}{|_{#2}}}
\newcommand{\T}[1]{{\mathcal{T}_{#1}}}
\newcommand{\Ttilde}[1]{{\widetilde{\mathcal{T}}_{#1}}}

\newcommand{\N}[1]{\mathcal{N}_#1}
\newcommand{\Ntilde}[1]{\tilde{\mathcal{N}}_#1}
\newcommand{\Hdiv}[1]{\bm{H}(\text{div},#1)}
\newcommand{\estim}[1]{\mathcal{E}_{#1}}
\newcommand{\estimf}[1]{\tilde{\mathcal{E}}_{#1}}
\newcommand{\Mha}{\bm{M}_h^{\bm{a}}}
\newcommand{\Qha}{Q_h^{\bm{a}}}
\newcommand{\wa}{{\omega_{\bm{a}}}}
\newcommand{\psia}{{\psi_{\bm{a}}}}
\newcommand{\lam}{\lambda_h^{\bm{a}}}

\newtheorem{prop}{Proposition}

\newtheorem{rem}{Remark}

\title{An equilibrated flux \emph{a posteriori} error estimator for defeaturing problems }

\author{Annalisa Buffa\thanks{MNS, Institute of Mathematics, École Polytechnique Fédérale de Lausanne, Switzerland \& IMATI CNR - Via Ferrata 5, 27100 Pavia  ({annalisa.buffa@epfl.ch})} \and Ondine Chanon\thanks{Institute of Analysis and Scientific Computing, Technische Universität Wien, Austria ({ondine.chanon@asc.tuwien.ac.at})} \and Denise Grappein\thanks{Dipartimento di Scienze Matematiche G. L. Lagrange, Politecnico di Torino, Italy, Member of GNCS INdAM Group ({denise.grappein@polito.it})} \and Rafael V\'azquez\thanks{Departamento de Matem\'atica Aplicada, Universidade de Santiago de Compostela, Spain ({rafael.vazquez@usc.es})} \and Martin Vohral\'ik \thanks{Inria, 2 rue Simone Iff, 75589 Paris, France \& CERMICS, Ecole des Ponts, 77455 Marne-la-Vallée, France ({martin.vohralik@inria.fr})}}

\begin{document}
	\maketitle
	\begin{abstract}
	An \textit{a posteriori} error estimator based on an equilibrated flux reconstruction is proposed for defeaturing problems in the context of finite element discretizations. Defeaturing consists in the simplification of a geometry by removing features that are considered not relevant for the approximation of the solution of a given PDE. In this work, the focus is on Poisson equation with Neumann boundary conditions on the feature boundary. The estimator accounts both for the so-called \textit{defeaturing error} and for the numerical error committed by approximating the solution on the defeatured domain. Unlike other estimators that were previously proposed for defeaturing problems, the use of the equilibrated flux reconstruction allows to obtain a sharp bound for the numerical component of the error. Furthermore, it does not require the evaluation of the normal trace of the numerical flux on the feature boundary: this makes the estimator well-suited for finite element discretizations, in which the normal trace of the numerical flux is typically discontinuous across elements.  The reliability of the estimator is proven and verified on several numerical examples. Its capability to identify the most relevant features is also shown, in anticipation of a future application to an adaptive strategy.
	\end{abstract}
		


	\section{Introduction}
The need of solving problems on complex domains, characterized by the presence of geometrical features of different scales and shapes, arises in many practical applications. In particular, in the process of simulation-based manufacturing, repeated simulations are often to be performed, in order to analyze the impact of design changes or to adjust geometric parameters. In many cases, before even solving the problem at hand, the first issue to overcome is the definition of the features themselves and the construction of a suitable computational mesh. For this reason it can be fundamental to simplify the geometry as much as possible, in order to avoid the definition of those features which may not have an actual impact on the accuracy of the solution. This process is commonly called \textit{defeaturing}. Some criteria based on some \textit{a priori} knowledge of the computational domain and of the properties of the materials have been used in the past (see, e.g., \cite{FINE00,FOUCAULT2004,THAKUR2009}). However, in order to fully automatize the process, an \textit{a posteriori} criterion is necessary and many different proposals can be found in literature (see \cite{CHOI2005,Ferrandes2009,GOP2007,GOP2009,LI2011b,LI2011,LI2013,LI2013b,SOK1999,TANG2013,TUR2009}).

In this paper we start from the work presented in \cite{BCV2022_arxiv,BCV2022}, which proposes an \emph{a posteriori} error estimator for analysis-aware defeaturing, in the context of the Poisson equation with Neumann boundary conditions on the feature boundary.  In particular in \cite{BCV2022_arxiv}, an estimator is designed to control the overall error between the exact solution of the PDE defined in the exact domain, and the numerical approximation of the solution of the corresponding
PDE defined in the defeatured domain. This estimator is made by two components, one accounting for the defeaturing error, i.e. the error committed by neglecting the features, and the other accounting for the numerical error committed when solving the problem on the defeatured geometry. The first component has the big advantage of being explicit with respect to the size of the geometrical features, and in \cite{BCV2022} the authors prove that it is a reliable and efficient bound for the energy norm of the defeaturing error. The second component is instead built as a residual-based estimator of the numerical error. The overall estimator is defined up to two positive parameters, related to the unknown constants appearing in the bounds of the defeaturing and of the numerical errors. Such parameters need to be tuned in order to correctly weight the two components. 

In order to partially overcome this issue, in this work we propose a novel \emph{a posteriori} error estimator that is strongly based on \cite{BCV2022} for what concerns its defeaturing component, but which resorts to an \textit{equilibrated flux reconstruction} (see, among others, \cite{AINSWORTH1997,DEST_MET1999,BraessScho2008,Luce_Wohlmuth}). Indeed, one of the main drawbacks of residual-based error estimators is that
the reliability constants are usually unknown and problem dependent. On the contrary, the difference between the numerical and the equilibrated flux provides an upper bound for the energy norm of the numerical error having reliability constant equal to 1. Although we do not get rid of the unknown constant related to the defeaturing component, the use of the equilibrated flux reconstruction also allows to avoid the computation of the normal trace of the numerical flux on the feature boundary. This makes the estimator well-suited for finite element discretizations, in which the normal trace of the numerical flux is typically not continuous.  On the contrary, the estimator proposed in \cite{BCV2022_arxiv} was designed to be applied along with an IGA discretization.

The equilibrated flux reconstruction is built following the steps in \cite{BraessScho2008,EV2015}, solving mixed local problems on patches of elements and leading to a discrete reconstructed flux in a Raviart---Thomas finite element space.

The paper is organized in four Sections. In Section~\ref{sec:not_and_mod} we introduce some notation and the defeaturing model problem, while in Section~\ref{sec:aposteriori_neg} and \ref{sec:aposteriori_pos} we derive and analyze an \emph{a posteriori} error estimator resorting to a generic equilibrated flux reconstruction and providing a bound for the overall error. Section \ref{sec:flux} describes a practical way to build the equilibrated flux reconstruction and, finally, in Section~\ref{sec:num_exp} the proposed estimator is validated by some numerical experiments.

	\section{Notation and model problem}\label{sec:not_and_mod}
	In the following we adopt the notation introduced in \cite{BCV2022}, which is here recalled for the sake of clarity. Let $\omega$ be any open $k$-dimensional manifold in $\mathbb{R}^d$, $d=2,3$ and $k\leq d$. We denote by $|\omega|$ the measure of $\omega$, and for any function $\varphi$ defined on $\omega$, we denote by $\overline{\varphi}^\omega$ its average over $\omega$. We will denote by $(\cdot,\cdot)_\omega$ the $L^2$-inner product on $\omega$ and by $||\cdot||_\omega$ the corresponding norm. If $k<d$, then $\langle \cdot,\cdot \rangle_\omega$ stands for a duality paring on $\omega$. For future use, let us define the quantity
	\begin{equation}
	c_\omega:=\begin{cases}
	\max(-\log(|\omega|),\zeta)^{\frac{1}{2}}&\text{if } k=1,~ d=2\\
	1 &\text{if } k=2, ~d=3
	\end{cases}\label{c_omega}
	\end{equation}
	where $\zeta \in \mathbb{R}$ is the unique solution of $\zeta=-\log(\zeta)$. 
	
	Let us consider an open Lipschitz domain $\Omega \subset \mathbb{R}^d$ and let us denote by $\partial \Omega$ its boundary. We suppose that $\Omega$ contains one feature $F\subset \mathbb{R}^d$, i.e. a geometrical detail of smaller scale, which is assumed to be an open Lipschitz domain, as well. The boundary of $F$ is denoted by $\partial F$. We consider two main types of features. In particular, a feature $F$ is said to be
	\begin{itemize}
		\item \textit{negative}, if $(\overline{F}\cap \overline{\Omega})\subset \partial \Omega$; 
		\item \textit{positive}, if $F\subset \Omega$.
	\end{itemize}
	In the following we will refer to $\Omega$ as the \textit{exact} or \textit{original} geometry.
	For the sake of simplicity we restrict ourselves to the case of an exact geometry with a single feature, but the generalization to the multiple feature case easily follows from \cite{AC_2023,BCV2022_arxiv}.
	
\begin{figure}
		\centering
		\begin{subfigure}[t]{.31\textwidth}
			\centering
			\includegraphics[width=1\linewidth]{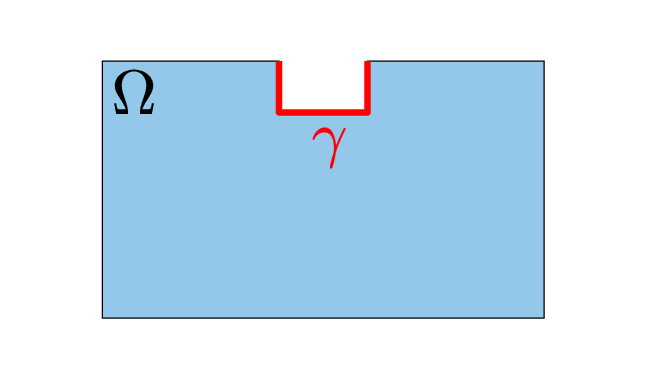}\hspace{-0.8cm}%
			\caption{Domain with negative feature.}
			\label{fig:ex_def_domain-a}
		\end{subfigure}\hspace{0.3cm}%
		\begin{subfigure}[t]{.31\textwidth}
			\centering
			\includegraphics[width=1\linewidth]{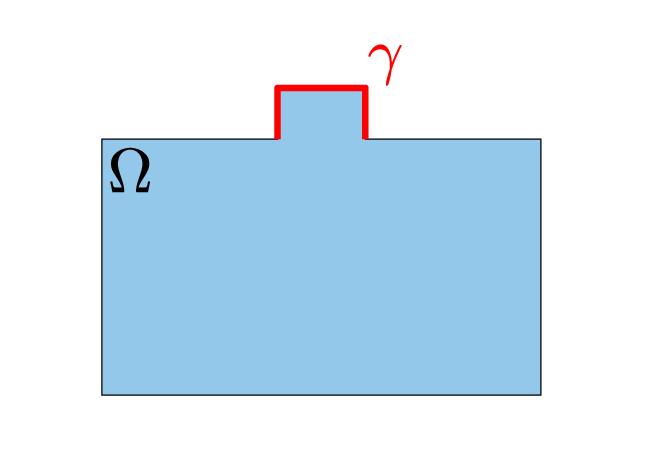}\hspace{-0.8cm}%
			\caption{Domain with positive feature.}
			\label{fig:ex_def_domain-b}
		\end{subfigure}\hspace{0.3cm}%
		\begin{subfigure}[t]{.31\textwidth}
			\centering
			\includegraphics[width=1\linewidth]{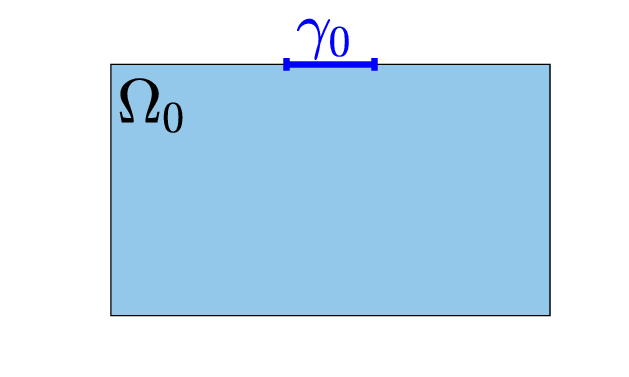}\hspace{-0.8cm}%
			\caption{Defeatured geometry.}
			\label{fig:ex_def_domain-c}
		\end{subfigure}%
		\caption{Domain with a negative feature, domain with a positive feature and corresponding defeatured domain for both configurations.}
		\label{fig:ex_def_domain}
	\end{figure}

	Let us now define the so called \textit{defeatured geometry}, i.e. 
	$\Omega_0\subset \mathbb{R}^d$ such that
	$$\Omega_0:=\begin{cases}
	\text{int}(\overline{\Omega}\cup \overline{F})&\text{ if } F \text{ is negative}\\
	\Omega\setminus \overline{F}&\text{ if } F \text{ is positive}.
	\end{cases}$$ Hence, if the feature is negative, $\Omega \subset \Omega_0$ (Figure~\ref{fig:ex_def_domain-a}), while if the feature is positive, $\Omega_0 \subset \Omega$ (Figure~\ref{fig:ex_def_domain-b}). In the following, the boundary of $\Omega_0$ is denoted by $\partial \Omega_0$.
	
	We denote by $\bm{n}$, $\bm{n_0}$ and $\bm{n_F}$ the unitary outward normals respectively of $\Omega$, $\Omega_0$ and $F$. Let $\partial \Omega=\overline{\Gamma_{\mathrm{D}}}\cup \overline{\Gamma_{\mathrm{N}}}$, with $\Gamma_{\mathrm{D}} \cap \Gamma_{\mathrm{N}}=\emptyset$ and $\Gamma_{\mathrm{D}} \neq \emptyset$, and we assume that $\partial F\cap \Gamma_{\mathrm{D}}=\emptyset$. Let $\gamma_0:=\partial F \setminus \overline{\Gamma_{\mathrm{N}}}\subset \partial \Omega_0$ and, finally, let $\gamma:=\partial F\setminus \overline{\gamma_0}\subset \partial \Omega$, so that $\partial F=\overline{\gamma_0}\cup\overline{\gamma}$ and $\gamma_0\cap \gamma=\emptyset$. Let us observe that, if $\gamma_0=\emptyset$, then we are in the case of a negative internal feature, i.e. $\Omega$ is a perforated domain (see, for an example, Figures~\ref{fig:meshT1} and \ref{fig:sol_ex_test3} in Section~\ref{sec:num_exp}).
	
	On the exact geometry $\Omega$ we use the Poisson problem as a model problem:
	\begin{equation}
	\begin{cases}
	-\Delta u =f & \text{in}~ \Omega\\
	u=g_{\mathrm{D}} &\text{on}~ \Gamma_{\mathrm{D}}\\
	\nabla u \cdot \bm{n}=g & \text{on}~ \Gamma_{\mathrm{N}},
	\end{cases} \label{eq:strong_omega}
	\end{equation}
	to which we will also refer as the \textit{original problem}.
	Defining
	$$H_{0,\Gamma_{\mathrm{D}}}^1(\Omega)=\big\lbrace v \in H^1(\Omega):~\tr{v}{\Gamma_{\mathrm{D}}}=0 \big\rbrace, \quad H_{g_{\mathrm{D}},\Gamma_{\mathrm{D}}}^1(\Omega)=\big\lbrace v \in H^1(\Omega):~\tr{v}{\Gamma_{\mathrm{D}}}=g_{\mathrm{D}} \big\rbrace, $$
	the variational formulation of Problem \eqref{eq:strong_omega} reads: \textit{find $u \in H^1_{g_{\mathrm{D}},\Gamma_{\mathrm{D}}}(\Omega)$ which satisfies}
	\begin{equation}
	(\nabla u,\nabla v)_{\Omega}=(f,v)_{\Omega}+\langle g,v\rangle_{\Gamma_{\mathrm{N}}} \quad \forall v \in H^1_{0,\Gamma_{\mathrm{D}}}(\Omega). \label{eq:weak_omega}
	\end{equation}
	On the defeatured geometry $\Omega_0$ we consider instead the problem
	\begin{equation}
	\begin{cases}
	-\Delta u_0 =f & \text{in}~ \Omega_0\\
	u_0=g_{\mathrm{D}} &\text{on}~ \Gamma_{\mathrm{D}}\\
	\nabla u_0 \cdot \bm{n}_0=g & \text{on}~ \Gamma_{\mathrm{N}}\setminus \gamma\\
	\nabla u_0 \cdot \bm{n}_0=g_0 & \text{on}~ \gamma_0
	\end{cases} \label{eq:strong_omega0}
	\end{equation} to which we will also refer to as \textit{defeatured problem}.
	With an abuse of notation, in the negative feature case, we denote by $f\in L^2(\Omega_0)$ a suitable $L^2$-extension of $f \in L^2(\Omega)$ to $F$, while the Neumann datum $g_0$ has to be chosen. The variational formulation of problem \eqref{eq:strong_omega0} reads: \textit{find $u_0 \in H^1_{g_{\mathrm{D}},\Gamma_{\mathrm{D}}}(\Omega_0)$ which satisfies, $\forall v \in H^1_{0,\Gamma_{\mathrm{D}}}(\Omega_0)$}
	\begin{equation}
	(\nabla u_0,\nabla v)_{\Omega_0}=(f,v)_{\Omega_0}+\langle g,v\rangle_{\Gamma_{\mathrm{N}}\setminus \gamma}+\langle g_0,v\rangle_{\gamma_0}. \label{eq:weak_omega0}
	\end{equation}
	
	Let us consider a partition $\T{h}$ of $\Omega_0$ consisting of closed triangles $K$ for $d=2$, or tetrahedrons for $d=3$, such that $\overline{\Omega_0}=\bigcup_{K \in \T{h}}K$. Hereby, we suppose that the mesh faces match with the boundaries $\Gamma_{\mathrm{D}}$, $\Gamma_{\mathrm{N}}\setminus \gamma$ and $\gamma_0$. Let us then introduce the set
	\begin{equation}Q_h=\mathcal{P}_p(\T{h}):=\left\lbrace q_h\in L^2(\Omega_0):~\tr{q_h}{K}\in \mathcal{P}_p(K),~\forall K \in \T{h} \right\rbrace,\label{eq:Qh}
	\end{equation} 
	with $\mathcal{P}_p(K)$ denoting the set of polynomials of degree at most $p\geq 1$ on $K \in \T{h}$, and
	$$V_h^0=\lbrace q_h \in\mathcal{C}^0(\overline{\Omega_0})\cap Q_h:~ \tr{q_h}{\Gamma_{\mathrm{D}}}=0\rbrace,\quad
	V_h:=\left\lbrace q_h\in \mathcal{C}^0(\overline{\Omega_0}):~\tr{q_h}{\Gamma_{\mathrm{D}}}=g_\mathrm{D}\right\rbrace.
	$$
	In the following, for the sake of simplicity, we assume $f\in Q_h$. Similarly, let us consider the partition of $\partial \Omega_0$ induced by the elements of $\T{h}$ and let us denote its restriction to $(\Gamma_{\mathrm{N}}\setminus \gamma)\cup\gamma_0$ by $\partial \Omega_{0,h}^\mathrm{N}$. Introducing
	$$g_{\mathrm{N}}=\begin{cases}
	g &\text{ on } \Gamma_{\mathrm{N}}\setminus \gamma\\
	g_0 &\text{ on } \gamma_0
	\end{cases} $$ we assume $g_{\mathrm{N}}$ to be an element of the broken space $\mathcal{P}_p(\partial \Omega_{0,h}^\mathrm{N})$, defined in the same manner as \eqref{eq:Qh}. Hence, the finite element approximation of \eqref{eq:weak_omega0} reads as:
	\textit{find $u_0^h \in V_h$ which satisfies, $\forall v_h \in V_h^0$}
	\begin{equation}
	(\nabla u_0^h,\nabla v_h)_{\Omega_0}=(f,v_h)_{\Omega_0}+\langle g,v_h\rangle_{\Gamma_{\mathrm{N}}\setminus \gamma}+\langle g_0,v_h\rangle_{\gamma_0}. \label{eq:discr_omega0}
	\end{equation}
	
	Let us remark that our aim is to never solve Problem~\eqref{eq:strong_omega}, but to design a proper \textit{a posteriori} error estimator capable to control the energy norm of the error committed by approximating the exact solution of \eqref{eq:strong_omega} by $u_0^h$. We will refer to this error as the \textit{overall error}, as it accounts both for the error introduced by defeaturing and for the error introduced by the numerical approximation of $u_0$. In particular, we aim at designing an estimator based on an \textit{equilibrated flux reconstruction}, which has the advantage of bounding the numerical error with a sharp reliability constant equal to 1. The flux reconstruction will be used to bound also the defeaturing error even if, in this case, we will not get rid of the unknown constant.  In the following we provide the definition of the overall error for the negative and positive feature case, referring again to \cite{BCV2022}.\\
	
	\textit{Negative feature:} in this case $\Omega\subset\Omega_0$, hence we restrict $u_0$ to $\Omega$ and we define the overall error as $||\nabla (u-\tr{u_0^h}{\Omega})||_\Omega$.\\
	
	\textit{Positive feature:} this case is slightly more complicated, since $u_0$ and its finite element approximation are defined only on $\Omega_0$ and $\Omega_0\subset \Omega$. Hence, in order to define the overall error, we need to extend $u_0$ to the feature $F$.
	\begin{figure}
		\centering
		\begin{subfigure}[t]{.45\textwidth}
			\centering
			\includegraphics[width=1\linewidth]{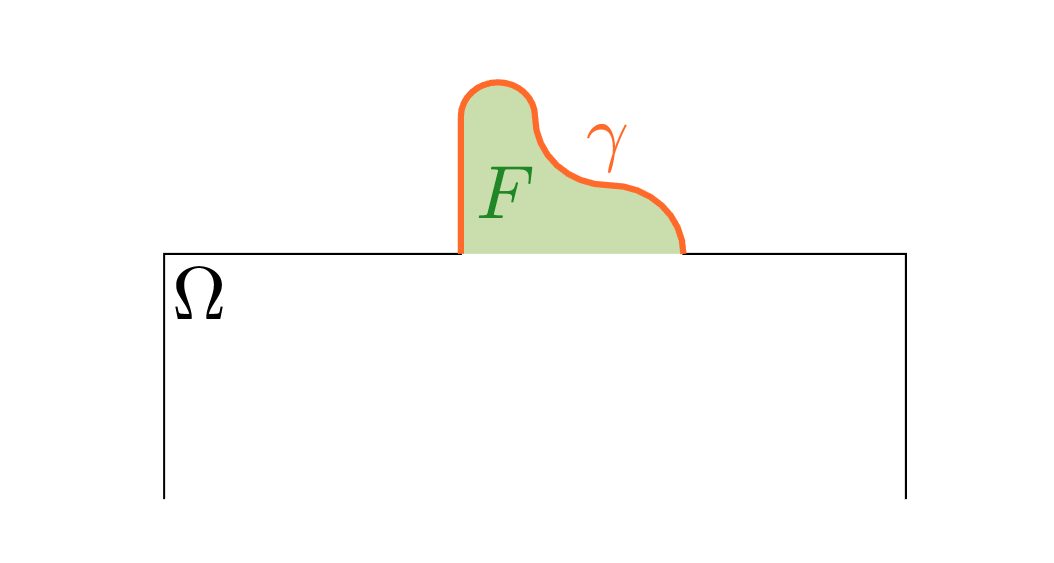}
			\caption{Domain with a positive feature $F$.}
			\label{F-a}
		\end{subfigure} \hspace{1cm}
		\begin{subfigure}[t]{.45\textwidth}
			\centering 
			\includegraphics[width=1\linewidth]{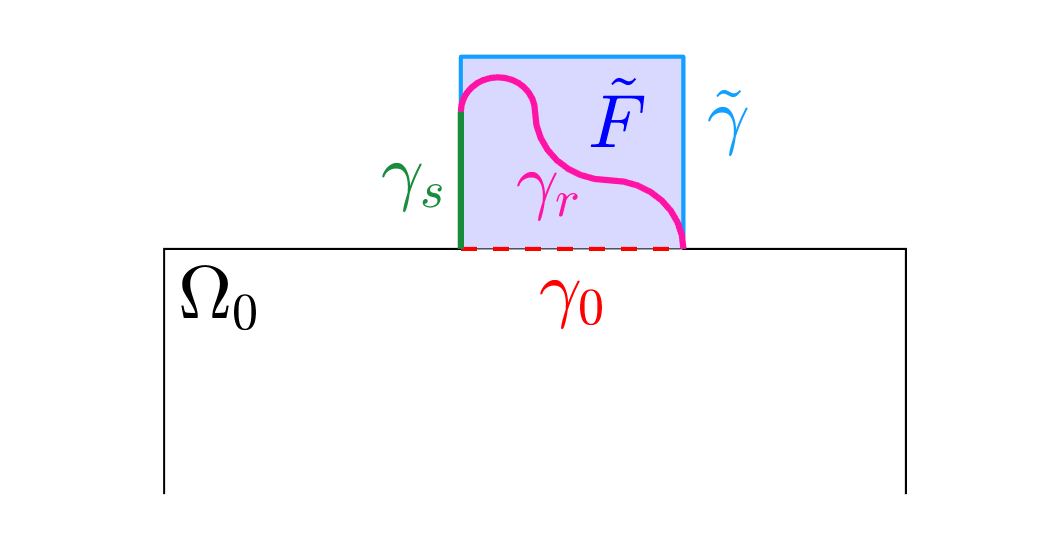}
			\caption{Extension $\tilde{F}$ of the positive feature and boundary nomenclature.}
			\label{F-b} %
		\end{subfigure}
		\caption{Example of geometry with a positive feature $F$ and extention of $F$ to $\tilde{F}$.}
		\label{fig:F}
	\end{figure}
	
	However meshing $F$ and solving a problem on it may be non trivial, in particular if $F$ has a complex boundary. Hence we follow the steps in \cite{BCV2022}: we consider a suitable extension $\tilde{F}\subset \mathbb{R}^d$ of $F$, being as simple as possible, in particular $F \subset \tilde{F}$ and $\gamma_0\subset(\partial \tilde{F} \cap \partial {F})$, as reported in Figure \ref{fig:F}. Let $\gamma$ be decomposed as $\gamma=\text{int}(\overline{\gamma_{\mathrm{s}}}\cup \overline{\gamma_{\mathrm{r}}})$ where $\gamma_{\mathrm{s}}=\gamma \cap \partial \tilde{F}$ is the portion of $\gamma$ shared by $\partial F$ and $\partial \tilde{F}$ and $\gamma_{\mathrm{r}}=\gamma\setminus \overline{\gamma_{\mathrm{s}}}$ (Figure \ref{F-b}). We denote by $\tilde{\bm{n}}$ the unitary outward normal of $\tilde{F}$ and we set $\tilde{\gamma}=\partial \tilde{F} \setminus \partial {F}$. On $\tilde{F}$ we solve the problem
	\begin{equation}
	\begin{cases}
	-\Delta \tilde{u}_0 =f & \text{in}~ \tilde{F}\\
	\tilde{u}_0=u_0 &\text{on}~ \gamma_0\\
	\nabla \tilde{u}_0 \cdot \tilde{\bm{n}}=\tilde{g} & \text{on}~ \tilde{\gamma}\\
	\nabla \tilde{u}_0 \cdot \tilde{\bm{n}}=g & \text{on}~ \gamma_{\mathrm{s}}
	\end{cases}\label{eq:strong_Ftilde}
	\end{equation}
	where, with an abuse of notation, we still denote by $f$ any $L^2$-extension of the forcing term to $\tilde{F}$ and the Neumann datum $\tilde{g}$ on $\tilde{\gamma}$ has to be chosen. Introducing
	$$H_{u_0,\gamma_0}^1(\tilde{F})=\left\lbrace v \in H^1(\tilde{F}):~\tr{v}{\gamma_0}=\tr{u_0}{\gamma_0} \right\rbrace, $$ the variational formulation of \eqref{eq:strong_Ftilde} is: \textit{find $\tilde{u}_0 \in H^1_{u_0,\gamma_0}(\tilde{F})$ which satisfies, $\forall v \in H^1_{0,\gamma_0}(\tilde{F})$}
	\begin{equation}
	(\nabla \tilde{u}_0,\nabla v)_{\tilde{F}}=(f,v)_{\tilde{F}}+\langle \tilde{g},v\rangle_{\tilde{\gamma}}+\langle g,v\rangle_{\gamma_{\mathrm{s}}}. \label{eq:weak_Ftilde}
	\end{equation}
	We denote by $\tilde{u}_0^h$ the finite element approximation of $\tilde{u}_0$ on a partition $\Ttilde{h}$ of $\tilde{F}$.  Note that this partition does not need to be conforming to $\gamma$. We suppose, however, that the mesh faces match with $\gamma_0$, $\tilde{\gamma}$ and $\gamma_{\mathrm{s}}$ and that $\T{h}$ matches with $\Ttilde{h}$ on $\gamma_0$. Let
	\begin{equation}
	\widetilde{Q}_h=\mathcal{P}_p(\Ttilde{h}):=\left\lbrace q_h\in L^2(\tilde{F}):~\tr{q_h}{K}\in \mathcal{P}_p(K),~\forall K \in \Ttilde{h}  \right\rbrace,\label{eq:Qhtilde}
	\end{equation} and let us introduce 
$$
	\widetilde{V}_h^0=\lbrace q_h \in\mathcal{C}^0(\overline{\tilde{F}})\cap \widetilde{Q}_h:~\tr{q_h}{\gamma_0}=0\rbrace, \quad
	\widetilde{V}_h:=\lbrace q_h\in \mathcal{C}^0(\overline{\tilde{F}})\cap \widetilde{Q}_h:~\tr{q_h}{\gamma_0}=\tr{u_0^h}{\gamma_0}\rbrace.
	$$
	We assume for simplicity that $\tr{f}{\tilde{F}}\in \widetilde{Q}_h$.
	Considering the partition of $\partial \tilde{F}$ induced by the elements of $\Ttilde{h}$ and denoting its restriction to $ \gamma_{\mathrm{s}}\cup\tilde{\gamma}$ as $\partial \tilde{F}_h^\mathrm{N}$, we also assume that
	$$\tilde{g}_\mathrm{N}=\begin{cases}
	g &\text{ on } \gamma_{\mathrm{s}}\\
	\tilde{g} &\text{ on } \tilde{\gamma}
	\end{cases} $$  is an element of the broken space $\mathcal{P}_p(\partial \tilde{F}_h^\mathrm{N})$.
	The finite element approximation of Problem \eqref{eq:weak_Ftilde} is hence: \textit{find $\tilde{u}_0^h \in \widetilde{V}_h$ which satisfies, $\forall v \in \widetilde{V}_h^0$}
	\begin{equation}
	(\nabla \tilde{u}_0^h,\nabla v_h)_{\tilde{F}}=(f,v_h)_{\tilde{F}}+\langle \tilde{g},v_h\rangle_{\tilde{\gamma}}+\langle g,v_h\rangle_{\gamma_{\mathrm{s}}}. \label{eq:discr_Ftilde}
	\end{equation}Finally, we define the extended defeatured solution and its numerical approximation as
	\begin{equation}
	u_{\mathrm{d}}:=\begin{cases}
	u_0&\text{ in }\Omega_0\\
	\tilde{u}_0&\text{ in } \tilde{F},
	\end{cases} \qquad u_{\mathrm{d}}^h:=\begin{cases}
	u_0^h&\text{ in }\T{h}\\
	\tilde{u}_0^h&\text{ in } \Ttilde{h},
	\end{cases}\label{eq:ud_udh}
	\end{equation}
	while the overall error is $||\nabla (u-u_{\mathrm{d}}^h)||_\Omega$.
	
	\section{Negative feature \emph{a posteriori} error estimator}\label{sec:aposteriori_neg}
	In this section we propose a reliable estimator for the overall error $||\nabla(u-\tr{u_0^h}{\Omega})||_{\Omega}$, in the case of a single negative feature. To simplify the notation, in the following we omit the explicit restriction of $u_0^h$ (and $u_0$) to $\Omega$. 
	
	Let us consider the solution to problem \eqref{eq:weak_omega0}: introducing the \textit{flux} $\bm{\sigma}=-\nabla u_0$, we have that $\bm{\sigma} \in \bm{H}(\text{div},\Omega_0)$, $\nabla \cdot \bm{\sigma}=f$, $-\bm{\sigma}\cdot \bm{n}_0=g$ on $\Gamma_{\mathrm{N}} \setminus \gamma$ and $-\bm{\sigma}\cdot \bm{n}_0=g_0$ on $\gamma_0$. At discrete level, a suitable definition of flux is more involved. Indeed $\nabla u_0^h\notin \Hdiv{\Omega_0}$ and hence the divergence equation and the Neumann boundary condition are not exactly satisfied. The idea behind the equilibrated flux reconstruction is to use the discrete solution $u_0^h$ to build a discrete flux $\flux$ such that $\flux \in \Hdiv{\Omega_0}$ is an approximation of $\bm{\sigma}$ satisfying
	\begin{equation}
	\begin{cases}
	\nabla \cdot \flux=f & \text{ in } \Omega_0\\
	\flux \cdot\bm{n}_0=-g & \text{ on } \Gamma_{\mathrm{N}} \setminus \gamma\\
	\flux \cdot\bm{n}_0=-g_0 & \text{ on }  \gamma_0.
	\end{cases}\label{eq:flux_prop}
	\end{equation}
	We will give more details about how an equilibrated flux reconstruction can actually be computed in Section~\ref{sec:flux}, following \cite{BraessScho2008,EV2015}. For the time being we assume that we have $\flux$ at our disposal.
	
	Referring to the notation introduced in Section~\ref{sec:not_and_mod}, let us introduce, on  $\gamma$, the quantity
	$$d_\gamma^h:=g+\flux\cdot \bm{n} \quad \text{on } \gamma$$
	which is the error between the Neumann datum $g$ on $\gamma$ and the normal trace of the equilibrated flux reconstruction. Following \cite{BCV2022}, denoting by $\overline{d_\gamma^h}^\gamma$ the average of $d_\gamma^h$ over $\gamma$, let us define
	\begin{equation}
	\estim{\gamma}:=\left(|\gamma|^\frac{1}{d-1}\left|\left|d_\gamma^h-\overline{d_\gamma^h}^\gamma\right|\right|_{\gamma}^2+c_\gamma^2|\gamma|^\frac{d}{d-1}\left|\overline{d_\gamma^h}^\gamma\right|^2\right)^{\frac{1}{2}}\label{eq:estim_gamma},
	\end{equation}
	\begin{equation}
	\estim{0}:=||\flux+\nabla u_0^h||_{\Omega_0}\label{eq:estim_num}
	\end{equation}
	where $c_\gamma$ is defined as in \eqref{c_omega}. Let us remark that, unlike  \cite{BCV2022}, the quantity $d_\gamma^h$ does not depend on the normal trace of the numerical flux, but on the normal trace of the equilibrated flux reconstruction, which is continuous across the elements of the mesh $\T{h}$. 
	The following proposition establishes our \emph{a posteriori} bound:
	\begin{prop}\label{prop:neg}
		Let $u$ be the solution of \eqref{eq:weak_omega} and $u_0^h$ the solution of \eqref{eq:discr_omega0}. Then
		\begin{equation}
		||\nabla(u-u_0^h)||_\Omega \leq C_D\estim{\gamma}+\estim{0},\label{eq:estim_neg_prop}
		\end{equation}
		with $C_D>0$ being a constant independent of the size of feature $F$.
	\end{prop}
	\begin{proof}
		Let $v \in H_{0,\Gamma_{\mathrm{D}}}^1(\Omega)$. Adding and subtracting $(\flux,\nabla v)_{\Omega}$, exploiting \eqref{eq:weak_omega}, applying Green's theorem and using the characterization of $\flux$ provided in \eqref{eq:flux_prop}, we have
		\begin{align}
		(\nabla(u-u_0^h),\nabla v)_{\Omega}&=(\nabla u+\flux,\nabla v)_{\Omega}-(\flux +\nabla u_0^h,\nabla v)_{\Omega}\nonumber\\
		&=(f-\nabla \cdot \flux,v)_{\Omega}+\langle g +\flux \cdot \bm{n},v\rangle_{\Gamma_{\mathrm{N}}}-(\flux +\nabla u_0^h,\nabla v)_{\Omega}
		\nonumber\\&=\langle g +\flux \cdot \bm{n},v\rangle_{\gamma}-(\flux +\nabla u_0^h,\nabla v)_{\Omega}\nonumber\\
		&=\langle d_\gamma^h,v\rangle_{\gamma}-(\flux +\nabla u_0^h,\nabla v)_{\Omega} \label{eq:def_error_preliminary}.
		\end{align}
		Referring the reader to the steps reported in \cite{BCV2022}, with the difference that the numerical flux is here substituted by the equilibrated flux reconstruction, it is possible to prove that 
		\begin{equation}
		\langle d_\gamma^h ,v\rangle_\gamma\leq
		C_D\estim{\gamma}||\nabla v ||_\Omega\label{eq:dmedio}
		\end{equation}
		with $C_D>0$ being a constant independent of the size of feature $F$ (see Theorem 4.3 in \cite{BCV2022}).
		If we choose $v=u-u_0^h\in H^1_{0,\Gamma_{\mathrm{D}}}(\Omega)$ in \eqref{eq:def_error_preliminary}, we apply \eqref{eq:dmedio} and the Cauchy--Schwarz inequality, we have
		\begin{align*}
		||\nabla(u-u_0^h)||_{\Omega}^2&\leq C_D \estim{\gamma}||\nabla(u-u_0^h)||_\Omega+||\flux +\nabla u_0^h||_\Omega||\nabla(u-u_0^h)||_\Omega\\&\leq C_D \estim{\gamma}||\nabla(u-u_0^h)||_\Omega+||\flux +\nabla u_0^h||_{\Omega_0}||\nabla(u-u_0^h)||_\Omega\\&= (C_D\estim{\gamma}+\estim{0})||\nabla(u-u_0^h)||_\Omega,
		\end{align*}
		where we have also exploited the fact that, in the negative feature case, $\Omega \subset \Omega_0$.
		Simplifying on both sides yields \eqref{eq:estim_neg_prop}.
	\end{proof}
	\begin{rem} It is well known from literature (see, among others, \cite{AINSWORTH1997,DEST_MET1999,BraessScho2008,Luce_Wohlmuth}) that the quantity $\estim{0}$ provides a sharp upper bound for the numerical error $||\nabla(u_0-u_0^h)||_{\Omega_0}$. Let us remark that, if no feature is present, the same result is provided also by \eqref{eq:estim_neg_prop}. Indeed, if $\gamma=\emptyset$, then $u=u_0$, $\Omega=\Omega_0$ and \eqref{eq:estim_neg_prop} reduces to $$||\nabla( u_0-u_0^h)||_{\Omega_0}\leq ||\flux+\nabla u_0^h||_{\Omega_0}.$$ 
		For this reason we will refer to $\estim{0}$ as the \textit{numerical component} of the estimator and to $\estim{\gamma}$ as the \textit{defeaturing} component.
	\end{rem}
	
	\section{Positive feature \emph{a posteriori} error estimator}\label{sec:aposteriori_pos}
	In this section we propose a reliable estimator for the overall error $||\nabla(u-u_{\mathrm{d}}^h)||_{\Omega}$, in the case of a single positive feature $F$. For the sake of generality, we consider the case in which $F$ is embedded in a smooth extension $\tilde{F}$, as detailed in Section~\ref{sec:not_and_mod}. 
	
	Let us introduce an equilibrated flux reconstruction on $\tilde{F}$, i.e. a discrete flux $\fluxf\in \Hdiv{\tilde{F}}$ built somehow from $\tilde{u}_0^h$ such that 
	\begin{equation}
	\begin{cases}
	\nabla \cdot \fluxf=f & \text{ in } \tilde{F}\\
	\fluxf \cdot\bm{\tilde{n}}=-\tilde{g} & \text{ on } \tilde{\gamma}\\
	\fluxf \cdot\bm{\tilde{n}}=-g & \text{ on } \gamma_{\mathrm{s}}.
	\end{cases}\label{eq:flux_prop_Ftilde}
	\end{equation}
	Again the details on the construction of this flux will be provided in Section \ref{sec:flux} and, for the time being, we assume we have $\fluxf$.
	In this case we define on $\gamma_0$ the quantity 
	$$d_{\gamma_0}^h:=\fluxf\cdot {\bm{n}}_F-g_0 \quad \text{on }\gamma_0$$
	which approximates the jump in the normal derivative of $u_{\mathrm{d}}$ on $\gamma_0$, while on $\gamma_{\mathrm{r}}$ we define
	$${d}_{\gamma_{\mathrm{r}}}^h:=\fluxf\cdot {\bm{n_F}}+g \quad \text{on }\gamma_{\mathrm{r}}$$
	which is the error between the Neumann datum $g$ on $\gamma_{\mathrm{r}}$ and the normal trace of the equilibrated flux reconstruction computed on $\tilde{F}$. Again we observe how, unlike \cite{BCV2022}, the normal trace of the numerical flux is not involved in the definition of these quantities. Denoting by $\overline{d_{\gamma_0}^h}^{\gamma_0}$ the average of $d_{\gamma_0}^h$ on $\gamma_0$ and by $\overline{{d}_{\gamma_{\mathrm{r}}}^h}^{\gamma_{\mathrm{r}}}$ the average of ${d}_{\gamma_{\mathrm{r}}}^h$ on $\gamma_{\mathrm{r}}$ and following \cite{BCV2022}, let us introduce
	\begin{equation}
	\estimf{\gamma_0}:=\Big(|\gamma_0|^\frac{1}{d-1}\left|\left|d_{\gamma_0}^h-\overline{d_{\gamma_0}^h}^{\gamma_0}\right|\right|_{{\gamma_0}}^2+c_{\gamma_0}^2|{\gamma_0}|^\frac{d}{d-1}\left|\overline{d_{\gamma_0}^h}^{\gamma_0}\right|^2\Big)^\frac{1}{2}\label{eq:estim_gamma0_pos}
	\end{equation}
	\begin{equation}
	\estimf{\gamma_{\mathrm{r}}}:=\Big(|\gamma_{\mathrm{r}}|^\frac{1}{d-1}\left|\left|d_{\gamma_{\mathrm{r}}}^h-\overline{d_{\gamma_{\mathrm{r}}}^h}^{\gamma_{\mathrm{r}}}\right|\right|_{\gamma_{\mathrm{r}}}^2+c_\gamma^2|\gamma_{\mathrm{r}}|^\frac{d}{d-1}\left|\overline{d_{\gamma_{\mathrm{r}}}^h}^{\gamma_{\mathrm{r}}}\right|^2\Big)^{\frac{1}{2}}\label{eq:estim_gammatilde_pos},
	\end{equation}
	where
	$c_{\gamma_0}$ and $c_{\gamma_{\mathrm{r}}}$ are defined as in \eqref{c_omega}. Let us also define
	\begin{equation}
	\estimf{0}:=||\fluxf+\nabla \tilde{u}_0^h||_{\tilde{F}},\label{eq:estim_num_Ftilde_pos}
	\end{equation} and let us recall that
	$
	\estim{0}=||\flux+\nabla u_0^h||_{\Omega_0},
	$
	where $\flux$ is an equilibrated flux reconstruction defined in $\Omega_0$ as in \eqref{eq:flux_prop}.
	\begin{prop}
		Let $u$ be the solution of \eqref{eq:weak_omega} and $u_{\mathrm{d}}^h$ be defined as in \eqref{eq:ud_udh}. Then
		\begin{equation}
		||\nabla(u-u_{\mathrm{d}}^h)||_\Omega \leq  \tilde{C}_D(\estimf{\gamma_0}^2+\estimf{\gamma_{\mathrm{r}}}^2)^\frac{1}{2}+(\estimf{0}^2+\estim{0}^2)^\frac{1}{2}
		\label{eq:estim_pos}
		\end{equation}
		with $\tilde{C}_D$ being a constant independent of the size of feature $F$.
	\end{prop}
	\begin{proof}
		Let us consider the restriction to $\Omega_0$ of the solution $u$ of problem \eqref{eq:strong_omega}, which satisfies
		\begin{equation}
		\begin{cases}
		-\Delta \tr{u}{\Omega_0} =f & \text{in}~ \Omega_0\\
		\tr{u}{\Omega_0} =g_{\mathrm{D}} &\text{on}~ \Gamma_{\mathrm{D}}\\
		\nabla \tr{u}{\Omega_0}  \cdot \bm{n}=g & \text{on}~ \Gamma_{\mathrm{N}}\setminus \gamma\\
		\nabla \tr{u}{\Omega_0}  \cdot \bm{n}_0=\nabla u  \cdot \bm{n}_0 & \text{on}~ \gamma_0.
		\end{cases} \label{eq:strong_omega_restrOmega_0}
		\end{equation}
		Omitting the explicit restriction on $u$ to $\Omega_0$ the variational formulation of problem \eqref{eq:strong_omega_restrOmega_0} reads: \textit{find $u \in H^1_{g_{\mathrm{D}},\Gamma_{\mathrm{D}}}(\Omega)$ which satisfies, $\forall v \in H^1_{0,\Gamma_{\mathrm{D}}}(\Omega)$}
		\begin{equation}
		(\nabla u,\nabla v)_{\Omega_0}=(f,v)_{\Omega_0}+\langle g,v\rangle_{\Gamma_{\mathrm{N}}\setminus \gamma}+\langle \nabla u \cdot \bm{n}_0,v\rangle_{\gamma_0}. \label{eq:weak_omega_restrOmega0}
		\end{equation}
		Let $v \in H_{0,\Gamma_{\mathrm{D}}}^1(\Omega_0)$. Adding and subtracting $(\flux,\nabla v)_{\Omega_0}$, exploiting  \eqref{eq:weak_omega_restrOmega0}, applying Green's theorem and using the characterization of $\flux$ provided in \eqref{eq:flux_prop}, we have
		\begin{align}
		(\nabla(u-u_0^h),\nabla v)_{\Omega_0}&=(\nabla u+\flux,\nabla v)_{\Omega_0}-(\flux +\nabla u_0^h,\nabla v)_{\Omega_0}\nonumber\\
		&=(f-\nabla \cdot \flux,v)_{\Omega_0}+\langle g +\flux \cdot \bm{n},v\rangle_{\Gamma_{\mathrm{N}}\setminus \gamma}\nonumber\\&\quad+\langle \nabla u\cdot \bm{n_0} +\flux \cdot \bm{n_0},v\rangle_{\gamma_0}-(\flux +\nabla u_0^h,\nabla v)_{\Omega_0}
		\nonumber\\&=\langle \nabla u\cdot \bm{n_0} -g_0,v\rangle_{\gamma_0}-(\flux +\nabla u_0^h,\nabla v)_{\Omega_0}. \label{eq:defpos_step1}
		\end{align}
		
		In order to obtain an actual error indicator we need to estimate the quantity $\langle \nabla u\cdot \bm{n_0} -g_0,v\rangle_{\gamma_0}$ on the right hand side, and for this reason we must consider the error committed on the feature as well. Hence, let us consider the restriction to the positive feature $F$ of the solution $u$ of \eqref{eq:strong_omega}, satisfying
		\begin{equation}
		\begin{cases}
		-\Delta \tr{u}{F} =f & \text{in}~ F\\
		\nabla \tr{u}{F}  \cdot \bm{n}_F=g & \text{on}~ \gamma\\
		\nabla \tr{u}{F}  \cdot \bm{n}_F=\nabla u  \cdot \bm{n}_F & \text{on}~ \gamma_0,
		\end{cases} \label{eq:strong_omega_restrF}
		\end{equation}
		so that, omitting the explicit restriction of $u$ to $F$, $u\in H^1(F)$ is one of the infinitely many solutions, defined up to a constant, of 
		\begin{equation}
		(\nabla u,\nabla v)_{F}=(f,v)_{F}+\langle g,v\rangle_{\gamma}+\langle \nabla u \cdot \bm{n}_F,v\rangle_{\gamma_0} \quad \forall v \in H^1(F). \label{eq:weak_omega_restrF}
		\end{equation}
		Let $v \in H^1(F)$. Adding and subtracting $(\fluxf,\nabla v)_{F}$, exploiting  \eqref{eq:weak_omega_restrF}, applying Green's theorem and using the characterization of $\fluxf$ provided in \eqref{eq:flux_prop_Ftilde}, we have
		\begin{align}
		(\nabla(u-\tilde{u}_0^h), &\nabla v)_{F}=(\nabla u+\fluxf,\nabla v)_{F}-(\fluxf +\nabla \tilde{u}_0^h,\nabla v)_{F}\nonumber\\
		&=(f-\nabla \cdot \fluxf,v)_{F}+\langle g +\fluxf \cdot \bm{n_F},v\rangle_{\gamma_{\mathrm{r}} \cup \gamma_{\mathrm{s}}}\nonumber\\&\quad +\langle \nabla u\cdot \bm{n_F} +\fluxf \cdot \bm{n_F},v\rangle_{\gamma_0}-(\fluxf +\nabla \tilde{u}_0^h,\nabla v)_{F}
		\nonumber\\&=\langle g +\fluxf \cdot \bm{n_F},v\rangle_{\gamma_{\mathrm{r}}}+\langle \nabla u\cdot \bm{n_F} +\fluxf \cdot \bm{n_F},v\rangle_{\gamma_0}-(\fluxf +\nabla \tilde{u}_0^h,\nabla v)_{F}.
		\label{eq:defpos_step2}
		\end{align} 
		Choosing $v \in H_{0,\Gamma_{\mathrm{D}}}(\Omega)$, observing that $\tr{v}{\Omega_0}\in H_{0,\Gamma_{\mathrm{D}}}^1(\Omega_0)$ and $\tr{v}{F}\in H^1(F)$, and summing \eqref{eq:defpos_step1} and \eqref{eq:defpos_step2}, we obtain
		\begin{align}
		(\nabla(u-u_{\mathrm{d}}^h), \nabla v)_{\Omega}&=(\nabla(u-u_0^h), \nabla v)_{\Omega_0}+(\nabla(u-\tilde{u}_0^h), \nabla v)_{F}\nonumber\\&=\langle g +\fluxf \cdot \bm{n_F},v\rangle_{\gamma_{\mathrm{r}}}+\langle\fluxf \cdot \bm{n_F}-g_0,v\rangle_{\gamma_0}\nonumber\\&\quad-(\flux +\nabla u_0^h,\nabla v)_{\Omega_0}-(\fluxf +\nabla \tilde{u}_0^h,\nabla v)_{F}\nonumber\\&=\langle d_{\gamma_{\mathrm{r}}}^h ,v\rangle_{\gamma_{\mathrm{r}}}+\langle d_{\gamma_0}^h ,v\rangle_{\gamma_0}-(\flux +\nabla u_0^h,\nabla v)_{\Omega_0}-(\fluxf +\nabla \tilde{u}_0^h,\nabla v)_{F}\label{eq:defpos_step3},
		\end{align} 
		where we have used that $\bm{n_F}=-\bm{n_0}$ on $\gamma_0$.
		
		For the terms involving $d_{\gamma_{\mathrm{r}}}^h$ and  $d_{\gamma_0}^h$ we proceed similarly to the negative feature case: referring the reader to Theorem 5.5 in \cite{BCV2022}, it is possible to prove that
		there exists a constant $\tilde{C}_D>0$ such that
		\begin{align}
		\langle d_{\gamma_0}^h ,v\rangle_{\gamma_0}+	\langle d_{\gamma_{\mathrm{r}}}^h ,&v\rangle_{\gamma_{\mathrm{r}}}\leq\tilde{C}_D\big(\estimf{\gamma_{\mathrm{r}}}^2+\estimf{\gamma_0}^2\big)^{\frac{1}{2}}||\nabla v||_{\Omega}.\label{eq:d_dmedio_pos}
		\end{align}	
		If we choose $v=u-u_{\mathrm{d}}^h\in H_{0,\Gamma_{\mathrm{D}}}(\Omega)$ in \eqref{eq:defpos_step3}, we use \eqref{eq:d_dmedio_pos} and the Cauchy--Schwarz inequality we obtain
		\begin{align*}
		||\nabla(u-u_{\mathrm{d}}^h)||_\Omega^2 &\leq \tilde{C}_D\big(\estimf{\gamma_{\mathrm{r}}}^2+\estimf{\gamma_0}^2\big)^{\frac{1}{2}}||\nabla(u-u_{\mathrm{d}}^h)||_\Omega+||\flux+\nabla u_0^h||_{\Omega_0}||\nabla(u-u_{\mathrm{d}}^h)||_{\Omega_0}\\&\quad+||\fluxf+\nabla \tilde{u}_0^h||_{F}||\nabla(u-u_{\mathrm{d}}^h)||_{F}\\
		&\leq \tilde{C}_D\big(\estimf{\gamma_{\mathrm{r}}}^2+\estimf{\gamma_0}^2\big)^{\frac{1}{2}}||\nabla(u-u_{\mathrm{d}}^h)||_\Omega\\
		&+\big(||\flux+\nabla u_0^h||_{\Omega_0}^2+||\fluxf+\nabla \tilde{u}_0^h||_{F}^2)^\frac{1}{2}(||\nabla(u-u_{\mathrm{d}}^h)||_{\Omega_0}^2+||\nabla(u-u_{\mathrm{d}}^h)||_{F}^2)^\frac{1}{2} \\
		&\leq \Big(\tilde{C}_D\big(\estimf{\gamma_{\mathrm{r}}}^2+\estimf{\gamma_0}^2\big)^{\frac{1}{2}}+(\estim{0}^2+\estimf{0}^2)^\frac{1}{2}\Big)||\nabla(u-u_{\mathrm{d}}^h)||_{\Omega},
		\end{align*}
		where, in the last step, we have exploited the fact that $F\subseteq \tilde{F}$.
		The thesis follows by simplifying on both sides.
	\end{proof}
	\begin{rem}\label{rem:F_or_Ftilde}
		If the feature $F$ is relatively simple, there is no need to use an extension, and problem \eqref{eq:weak_Ftilde} is solved directly in $F$. In this case, maintaining the tilde notation (since $F=\tilde{F}$), expression \eqref{eq:estim_pos} simplifies into 	
		\begin{equation}
		||\nabla(u-u_{\mathrm{d}}^h)||_\Omega \leq \tilde{C}_D\estimf{\gamma_0}+(\estimf{0}^2+\estim{0}^2)^\frac{1}{2}.
		\end{equation}
	\end{rem}

	\section{Equilibrated fluxes reconstruction}\label{sec:flux}
	In this section we describe how to build, in practice, an equilibrated flux starting from the discrete solution of the defeatured problem $u_0^h$ or from $u_{\mathrm{d}}^h$ in the positive feature case. The proposed procedure is directly adapted from \cite{BraessScho2008,EV2015} and resorts to a local reconstruction of the fluxes.
	Given the triangular/tetrahedral mesh $\T{h}$ built on the defeatured geometry $\Omega_0$, let us denote by $\N{h}$ the set of its vertices and let us divide it into interior vertices $\N{h}^{\mathrm{int}}$ and boundary vertices $\N{h}^{\mathrm{ext}}$.
	We aim at reconstructing the flux in the Raviart--Thomas space of order $p\geq 1$, namely in
	$$\bm{M}_h:=\left\lbrace  \bm{v}_h\in \Hdiv{\Omega_0}:~\tr{\bm{v}_h}{K} \in [\mathcal{P}_p(K)]^d+\bm{x}\mathcal{P}_p(K),~\forall K \in \T{h}\right\rbrace.$$ 
	The best choice for the equilibrated flux reconstruction would then be 
	\begin{equation}
	\flux=\arg\min_{\bm{v}_h \in \bm{M}_h}||\bm{v}_h+\nabla u_0^h||_{\Omega_0} \quad \text{ subject to } \begin{cases}\nabla \cdot \bm{v}_h =f &\text{ in }\Omega_0\\
	\bm{v}_h \cdot \bm{n}=-g &\text{ on }\Gamma_{\mathrm{N}}\setminus \gamma\\
	\bm{v}_h \cdot \bm{n}_0=-g_0 &\text{on }\gamma_0.\\
	\end{cases}\label{eq:flux_global}
	\end{equation}
	However, finding $\flux$ through \eqref{eq:flux_global} implies solving a global optimization problem in the domain $\Omega_0$. 
	
	Following \cite{BraessScho2008,EV2015} we adopt instead a different strategy, in which local flux reconstructions are built on patches $\wa$ of elements sharing a vertex $\bm{a} \in \N{h}$. Let us denote by $\psi_{\bm{a}}$ the \textit{hat} function in $\mathcal{P}_1(\T{h})\cap H^1(\Omega_0)$ taking value 1 in vertex $\bm{a}$ and 0 on all the other vertices.  Let us denote by $\partial \wa$ the boundary of the patch $\wa$ and let $\partial \omega_{\bm{a}}^0\subseteq \partial \wa$ be defined as
	$$\partial \omega_{\bm{a}}^0=\lbrace \bm{x} \in \partial \wa~:\psi_{\bm{a}}(\bm{x})=0\rbrace,$$ and $\partial \omega_{\bm{a}}^\psi=\partial \wa \setminus {\partial \omega_{\bm{a}}^0}$. Let us remark that, if $\bm{a}\in \N{h}^{\mathrm{int}}$ then $\partial \omega_{\bm{a}}^0=\partial \wa$. Let $\Gamma_{\mathrm{N}}^0=(\Gamma_{\mathrm{N}}\setminus \gamma)\cup \gamma_0$ and let us introduce
	$$\bm{M}_h^{\bm{a},0}=\left\lbrace \bm{v}_h \in \bm{M}_h(\wa):~\bm{v}_h\cdot \bm{n}_{\wa}=0 \text{ on }\partial \omega_{\bm{a}}^0\cup \big(\partial \omega_{\bm{a}}^\psi\cap\Gamma_{\mathrm{N}}^0\big)\right\rbrace$$
	and
	\begin{align*}
	\Mha&:=\begin{cases}\bm{M}_h^{\bm{a},0}& \text{ if } \bm{a} \in \N{h}^{\mathrm{int}}\\[0.3em]
	\big\lbrace \bm{v}_h \in \bm{M}_h(\wa):~\bm{v}_h\cdot \bm{n}_{\wa}=0 \text{ on }\partial \omega_{\bm{a}}^0,\\\hspace{2.7cm} \bm{v}_h\cdot \bm{n}_{\wa}=-\psia g \text{ on }\partial\omega_{\bm{a}}^\psi\cap (\Gamma_{\mathrm{N}}\setminus \gamma),\\\hspace{2.7cm}\bm{v}_h\cdot \bm{n}_{\wa}=-\psia g_0 \text{ on } \partial\omega_{\bm{a}}^\psi\cap \gamma_0~\big\rbrace& \text{ if } \bm{a} \in \N{h}^{\mathrm{ext}},
	\end{cases}\\[1em]
	\Qha&:=\begin{cases}\big\lbrace q_h \in Q_h(\wa):(q_h,1)_{\wa}=0 \big\rbrace & \text{ if } \bm{a} \in \N{h}^{\mathrm{int}} \text{ or } \bm{a} \in \text{int}\big(\overline{\Gamma_{\mathrm{N}}^0}) \\[0.4em] Q_h(\wa) &\text{ if } \bm{a} \in \N{h}^{\mathrm{ext}} \text{ and } \bm{a} \notin  \text{int}\big(\overline{\Gamma_{\mathrm{N}}^0}\big),
	\end{cases}
	\end{align*}
	where  $\bm{M}_h(\wa)$ and $Q_h(\wa)$ are respectively the restrictions of $\bm{M}_h$ and $Q_h$ to the patch $\wa$ and $Q_h$ is defined as in \eqref{eq:Qh}.
	We then look for local equilibrated flux reconstructions as
	\begin{equation}
	\fluxa=\arg\min_{\bm{v}_h \in \Mha}||\bm{v}_h+\psi_{\bm{a}}\nabla u_0^h||_{0,\wa}\quad \text{ subject to }\quad  \nabla \cdot \bm{v}_h =\psia f-\nabla\psia\cdot \nabla u_0^h \label{eq:local_argmin}
	\end{equation} and then we set $\flux =\sum_{\bm{a}\in \N{h}}\fluxa.$
	
	The optimization problem \eqref{eq:local_argmin} is equivalent to look for $\fluxa \in \Mha$ and $\lam \in \Qha$ such that,
	\begin{align}
	(\fluxa,\bm{v}_h)_\wa-(\lam,\nabla \cdot \bm{v}_h)_\wa=-(\psia\nabla u_0^h,\bm{v}_h)_\wa \quad \forall \bm{v}_h \in \bm{M}_h^{\bm{a},0} \label{eq:mixed1}\\[0.5em]
	(\nabla \cdot \fluxa,q_h)_\wa=(\psia f,q_h)_\wa-(\nabla\psia\cdot \nabla u_0^h,q_h)_\wa \quad \forall q_h \in \Qha\label{eq:mixed2},
	\end{align}
	which is the strategy that we are actually adopting in practice.
	
	The equilibrated flux on the extension $\tilde{F}$ of a positive feature $F$ is reconstructed exactly in the same manner. Denoting by $\Ntilde{h}^{\mathrm{int}}$ and $\Ntilde{h}^{\mathrm{ext}}$ respectively the internal and the boundary vertices of the mesh $\Ttilde{h}$ of $\tilde{F}$, introducing
	\begin{align*}
	&\widetilde{\bm{M}}_h(\Ttilde{h}):=\big\lbrace  \bm{v}_h\in \Hdiv{\tilde{F}}:~\tr{\bm{v}_h}{K} \in [\mathcal{P}_p(K)]^d+\bm{x}\mathcal{P}_p(K),~\forall K \in \Ttilde{h}\big\rbrace,
	\end{align*} and recalling the definition of $\widetilde{Q}_h$ given in \eqref{eq:Qhtilde},
	we look for the couple $({\fluxf}^{\bm{a}},\tilde{\lambda}_h^{\bm{a}})$ in the sets and spaces
	\begin{align*}
	\widetilde{\bm{M}}_h^{\bm{a}}&:=\begin{cases}\left\lbrace \bm{v}_h \in \widetilde{\bm{M}}_h(\wa):~\bm{v}_h\cdot \bm{n}_{\wa}=0 \text{ on }\partial \omega_{\bm{a}}^0\cup \big(\partial \omega_{\bm{a}}^\psi\cap(\tilde{\gamma}\cup\gamma_s)\big)\right\rbrace& \text{ if } \bm{a} \in \tilde{\mathcal{N}}_{h}^{\mathrm{int}}\\[0.5em]
	\big\lbrace \bm{v}_h \in \widetilde{\bm{M}}_h(\wa):~\bm{v}_h\cdot \bm{n}_{\wa}=0 \text{ on } \partial \omega_{\bm{a}}^0,\\\hspace{2.7cm}\bm{v}_h\cdot \bm{n}_{\wa}=-\psia \tilde{g} \text{ on }\partial \omega_{\bm{a}}^\psi\cap \tilde{\gamma},\\
	\hspace{2.7cm}\bm{v}_h\cdot \bm{n}_{\wa}=-\psia g \text{ on }\partial \omega_{\bm{a}}^\psi\cap \gamma_{\mathrm{s}}~\big\rbrace& \text{ if } \bm{a} \in \tilde{\mathcal{N}}_{h}^{\mathrm{ext}},
	\end{cases}\\[0.8em]
	\widetilde{Q}_h^{\bm{a}}&:=\begin{cases}\big\lbrace q_h \in \widetilde{Q}_h(\wa):~(q_h,1)_{\wa}=0 \big\rbrace  &\text{ if } \bm{a} \in \tilde{\mathcal{N}}^{\mathrm{int}}_h \text{ or } \bm{a} \in \text{int}\big(\overline{\tilde{\gamma}\cup\gamma_\mathrm{s}}\big)\\[0.5em] \widetilde{Q}_h(\wa)&\text{ if } \bm{a} \in \tilde{\mathcal{N}}_{h}^{\mathrm{ext}} \text{ and }\bm{a} \notin \text{int}\big(\overline{\tilde{\gamma}\cup\gamma_\mathrm{s}}\big)
	\end{cases}
	\end{align*}
	solving a problem analogous to \eqref{eq:mixed1}-\eqref{eq:mixed2}.

	\section{Numerical experiments}\label{sec:num_exp}
	In this section we propose some numerical experiments to validate the proposed estimator. We here focus on the case $d=2$ and $p=1$. All the simulations were performed in Matlab and meshes were built using the Triangle mesh generator \cite{shewchuk96b}. For each element $K$ of a mesh $\T{h}$ we denote by $h_K$ the diameter of the element, and we choose as a mesh parameter $h=\max_{K \in \T{h}}h_K$.
	
	In the following we define the \textit{effectivity index} as the ratio between the total estimator and the overall error. i.e.
	\begin{equation}
	\eta=\frac{\estim{\mathrm{tot}}}{||\nabla (u-u_0^h)||_{\Omega}} \label{eta}.
	\end{equation} 
	In case of a domain $\Omega$ characterized by a single negative feature, following Proposition~\ref{prop:neg}, the total estimator is defined as 
	\begin{equation}
	\estim{\mathrm{tot}}=C_D\estim{\gamma}+\estim{0}, \label{totestim_neg}
	\end{equation}
	while in the case of a single positive feature it is
	\begin{equation}
	\estim{\mathrm{tot}}=\tilde{C}_D\estimf{\gamma_0}+(\estimf{0}^2+\estim{0}^2)^\frac{1}{2}, \label{totestim_pos_F} 
	\end{equation}
	assuming $\tilde{F}=F$ (see Remark \ref{rem:F_or_Ftilde}). 
	For all the proposed experiments, a reference solution is built by solving the problem on the original geometry $\Omega$ by linear finite elements on a very fine mesh. With an abuse of notation this reference solution is still denoted by $u$.
	
	Three numerical examples are proposed. In Test 1 we consider the case of a single negative internal feature, analyzing the convergence of the total estimator and of the overall error under mesh refinement and feature size reduction. Test 2 deals instead with the case of positive and negative boundary features and the convergence of the estimator and of the error are again analyzed. Finally, in Test 3 we consider the presence of multiple internal negative features, showing how the proposed estimator allows to point out which features have the greatest impact on the error.
	
	\subsection{Test 1: negative internal feature}
	\begin{figure}
		\centering
		\begin{subfigure}[t]{0.4\textwidth}
			\centering
			\includegraphics[width=0.95\linewidth]{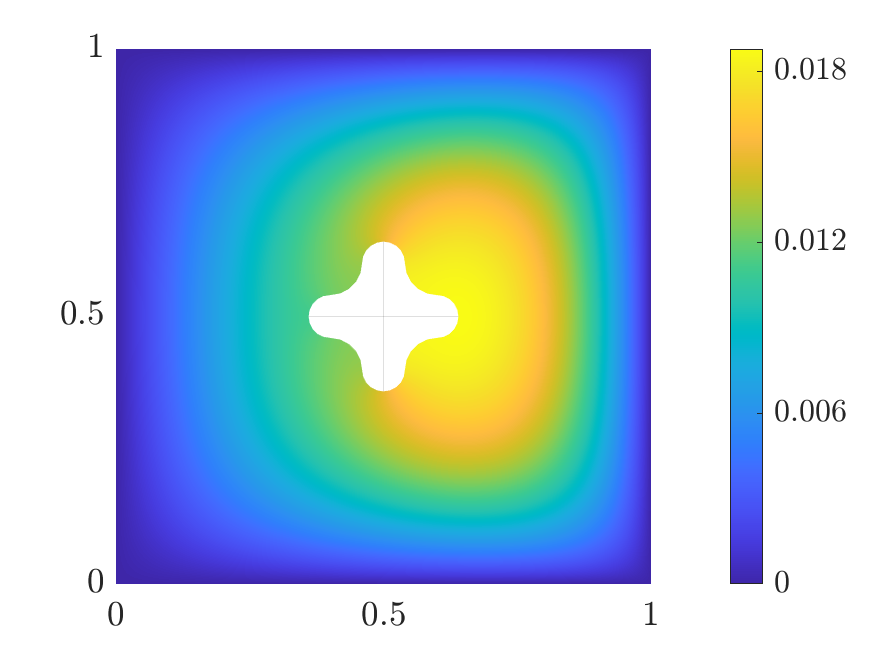}
			\caption{Reference solution for $\epsilon=0.14$.}
			\label{fig:meshT1-a}
		\end{subfigure}\hspace{0.7cm}%
		\begin{subfigure}[t]{0.4\textwidth}
			\centering
			\includegraphics[width=0.95\linewidth]{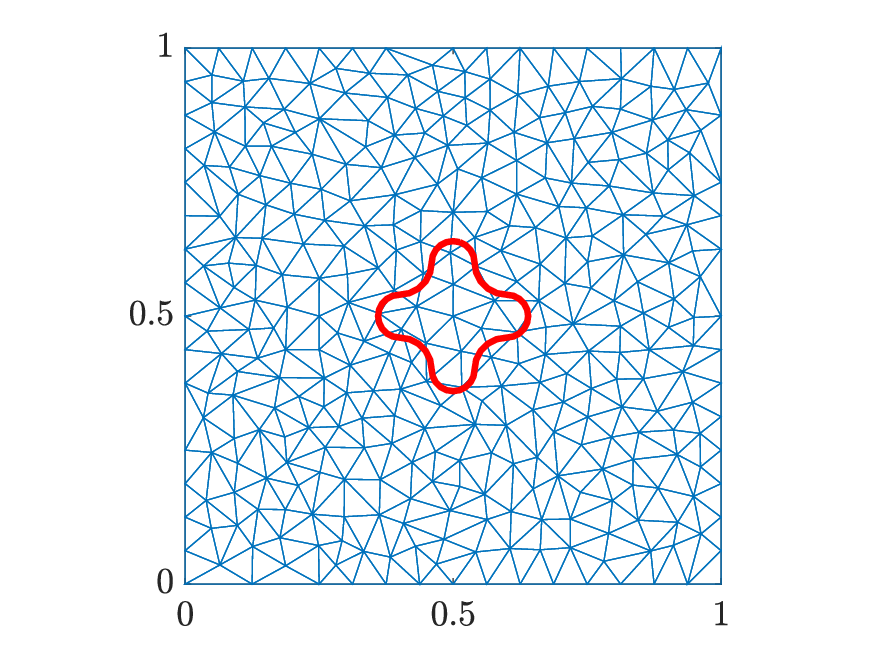}
			\caption{Example of computational mesh defined on the defeatured geometry.}
			\label{fig:meshT1-b}
		\end{subfigure}
		\caption{Test 1: reference solution $u_\epsilon$ ($\epsilon=0.14$), and defeatured geometry with an example of computational mesh used to compute the defeatured solution and the equilibrated flux reconstruction.}
		\label{fig:meshT1}
	\end{figure}
	
	For this first numerical example we consider a square domain characterized by a single negative (see Figure \ref{fig:meshT1-a}). We denote by $\epsilon$ the characteristic size of the feature, i.e. the radius of the circle circumscribing the feature itself.
	
	Setting $\Omega_0=(0,1)^2$ and $\Omega_\epsilon=\Omega_0\setminus \overline{F_\epsilon}$ we consider on the exact geometry $\Omega_\epsilon$ the problem
	\begin{equation}
	\begin{cases}
	-\Delta u_\epsilon=f &\text{in } \Omega_\epsilon\\
	u_\epsilon=0 &\text{on } \partial \Omega_\epsilon \setminus \gamma_\epsilon\\
	\nabla u_\epsilon\cdot \bm{n}=0 &\text{on } \gamma_\epsilon,
	\end{cases}\label{prob2}
	\end{equation}
	with $f(x,y)=x$.
	For what concerns the defeatured problem, $\gamma_0=\emptyset$, since the feature is internal.
	
	Figure \ref{fig:meshT1-b} reports an example of the computational mesh $\T{h}$ used to solve the defeatured problem on $\Omega_0$. We remark how the mesh has no need to be conforming to the feature boundary since the equilibrated flux $\flux$ is reconstructed on the defeatured geometry, which is blind to the feature, and $\estim{\gamma}$ is computed by simply defining a proper quadrature rule on the feature boundary itself and evaluating the normal trace of $\flux$ in the chosen quadrature nodes.

\begin{figure}
		\centering
		\begin{subfigure}{0.315\textwidth}
			\centering
			\includegraphics[width=1.15\linewidth]{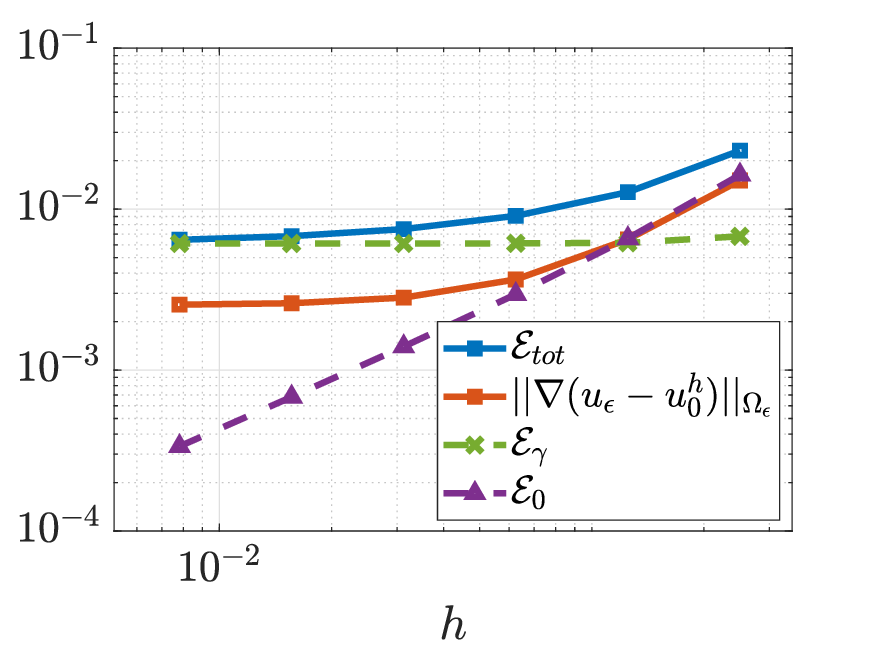}
			\caption{$\epsilon=7.00e-2$}
			\label{fig:varh_test1-a}
		\end{subfigure}
		\hfill
		\begin{subfigure}{0.315\textwidth}
			\centering
			\includegraphics[width=1.15\linewidth]{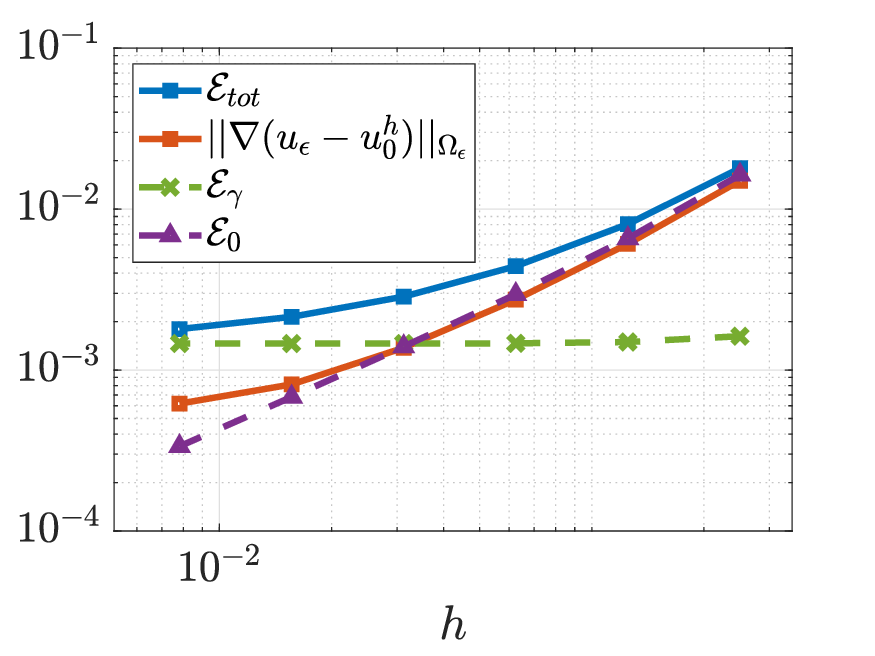}
			\caption{$\epsilon=1.75e-2$}
			\label{fig:varh_test1-b}
		\end{subfigure}
		\hfill
		\begin{subfigure}{0.315\textwidth}
			\centering
			\includegraphics[width=1.15\linewidth]{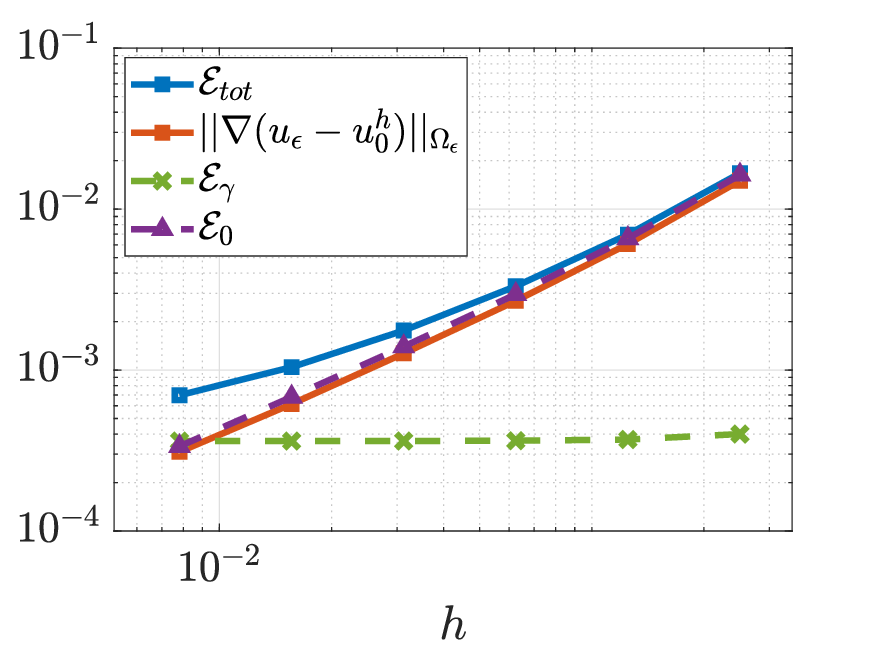}
			\caption{$\epsilon=4.38e-3$}
			\label{fig:varh_test1-c}
		\end{subfigure}
		\hfill
		\caption{Test 1: energy norm of the error (full red line), total estimator (full blue line) and its components (dashed lines) under mesh refinement and for three fixed feature sizes}
		\label{fig:varh_test1}
	\end{figure}

	In the numerical experiments that follow we choose $C_D=1$ in \eqref{totestim_neg}.
	Figure~\ref{fig:varh_test1} shows the convergence of the estimator $\estim{\mathrm{tot}}$ and of the energy norm of the overall error $||\nabla(u_\epsilon-u_0^h)||_{\Omega}$, under mesh refinement. The values of $\estim{\gamma}$ and $\estim{0}$ are also reported. Three fixed values of $\epsilon$ are considered, namely $\epsilon=7.00\cdot 10^{-2}, ~1.75\cdot 10^{-2},~4.83\cdot 10^{-3}$. As expected, going from a coarse to a fine mesh, the error reaches a plateau when the defeaturing error becomes more relevant than the numerical one. The bigger the feature is, the earlier the plateau is reached. The same behavior is captured also by the estimator. For a fixed feature size, the value of $\estim{\gamma}$ remains constant while $\estim{0}$ converges as $\mathcal{O}(h)$, as expected.
	
\begin{figure}
		\centering
		\begin{subfigure}[t]{0.4\textwidth}
			\centering
			\includegraphics[width=1.04\linewidth]{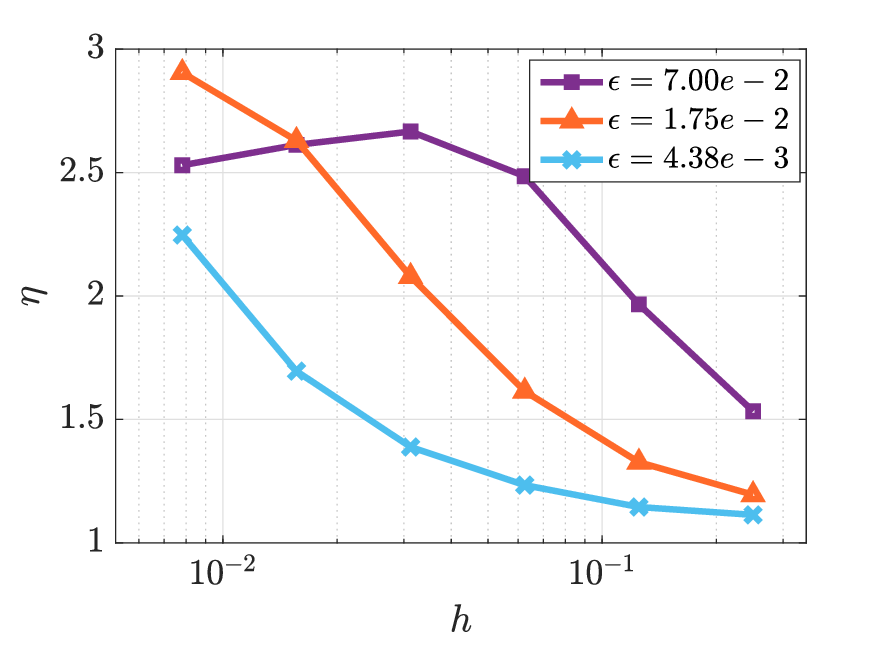}
			\caption{Effectivity index under mesh refinement, for three fixed feature sizes.}
			\label{fig:eff_test2-a}
		\end{subfigure}\hspace{0.7cm}%
		\begin{subfigure}[t]{0.4\textwidth}
			\centering
			\includegraphics[width=1\linewidth]{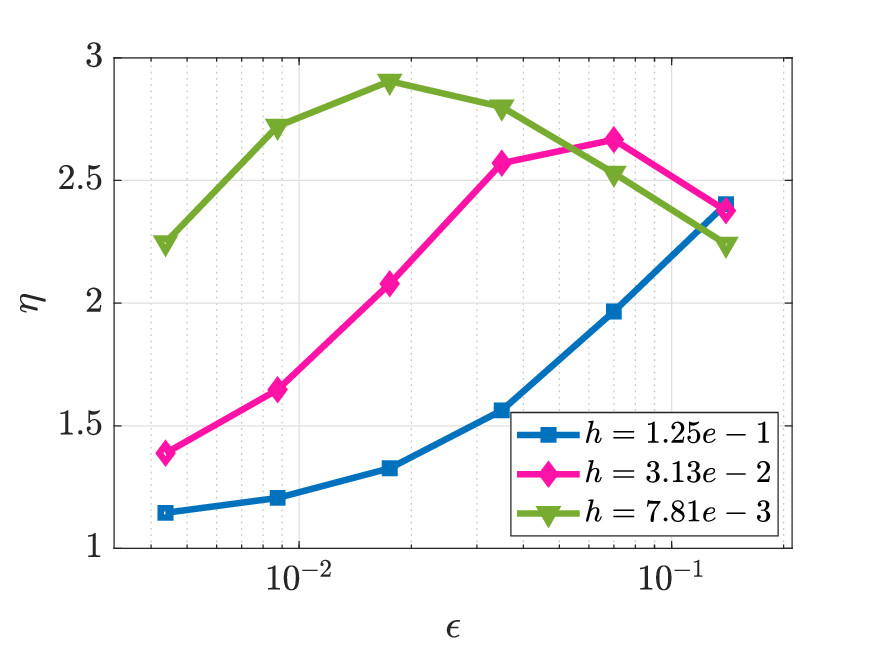}
			\caption{Effectivity index under reduction of feature size, for three fixed mesh sizes.}
			\label{fig:eff_test2-b}
		\end{subfigure}
		\caption{Test 1: Effectivity index under mesh refinement and under feature size reduction.}
		\label{fig:eff_test2}
	\end{figure}

	The trend of the effectivity index $\eta$ under mesh refinement and for the three considered feature sizes is reported in Figure~\ref{fig:eff_test2-a}. As expected, since the numerical source of the error is sharply bounded by $\estim{0}$, when $\estim{0}\gg\estim{\gamma}$ we have $\eta \sim1$. The effectivity index appears instead to be around $2.5$ when the defeaturing component is dominating. The highest values of $\eta$, namely $\eta \sim 3$, are registered when $\estim{\gamma}>\estim{0}$ but the two components have still a comparable magnitude.
	
\begin{figure}
		\centering
		\begin{subfigure}{0.315\textwidth}
			\includegraphics[width=1.15\linewidth]{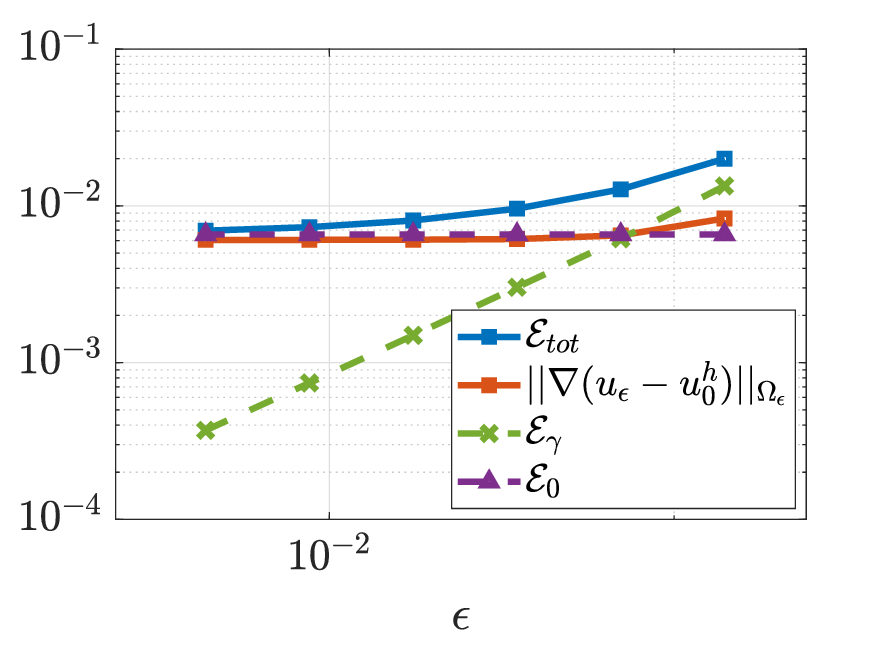}
			\caption{$h=1.25e-1$}
			\label{fig:vareps_test1-a}
		\end{subfigure}
		\hfill
		\begin{subfigure}{0.315\textwidth}
			\includegraphics[width=1.15\linewidth]{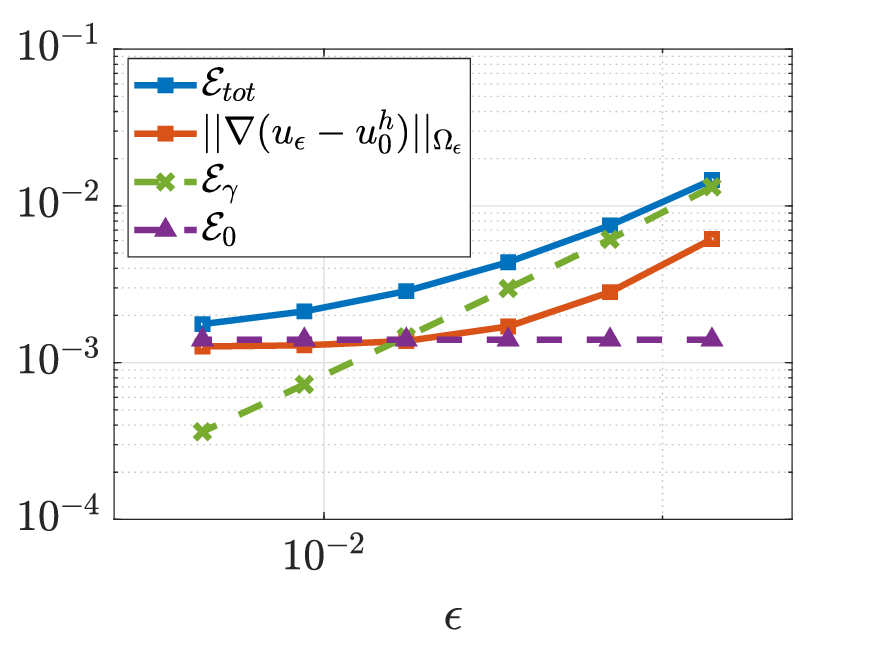}
			\caption{$h=3.13e-2$}
			\label{fig:vareps_test1-b}
		\end{subfigure}
		\hfill
		\begin{subfigure}{0.315\textwidth}
			\includegraphics[width=1.15\linewidth]{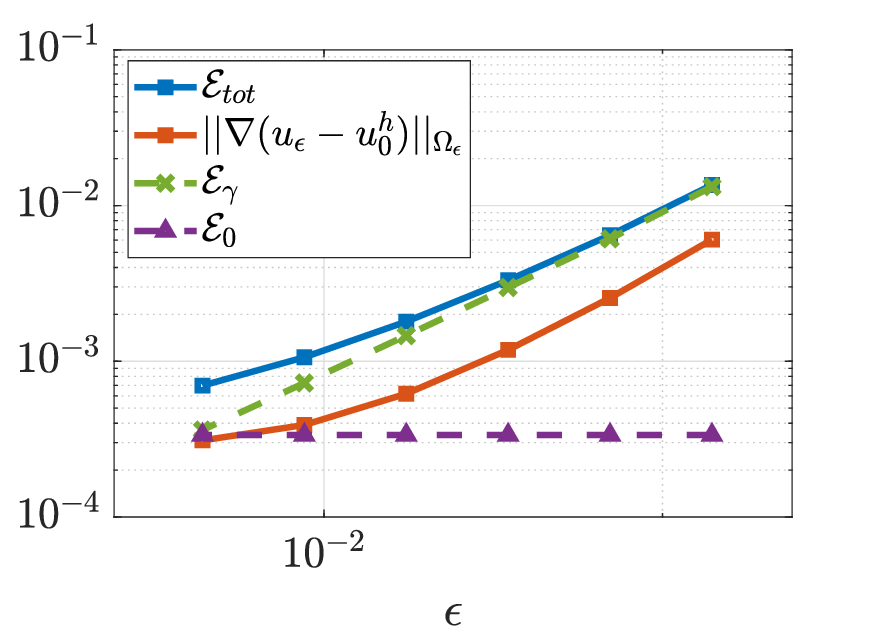}
			\caption{$h=7.81e-3$}
			\label{fig:vareps_test1-c}
		\end{subfigure}
		\hfill
		\caption{Test 1: energy norm of the error, total estimator and its components under the reduction of the feature size and for three different fixed mesh sizes.}
		\label{fig:vareps_test1}
	\end{figure}

	Figure~\ref{fig:vareps_test1} focuses instead on the convergence of the estimator and of the error under the reduction of the feature size, for three fixed mesh sizes, namely $h=1.25\cdot 10^{-1},~3.13\cdot 10^{-2},~7.81\cdot 10^{-3}$. As expected, both the error and the estimator reach a plateau when the numerical error dominates over the defeaturing one. The value of the effectivity index $\eta$ is reported in Figure~\ref{fig:eff_test2-b}, with the same considerations done for the Figure~\ref{fig:eff_test2-a} still holding.


	\subsection{Test 2: positive and negative boundary features}
	\begin{figure}
		\centering
		\begin{subfigure}[t]{0.4\textwidth}
			\centering
			\includegraphics[width=0.96\linewidth]{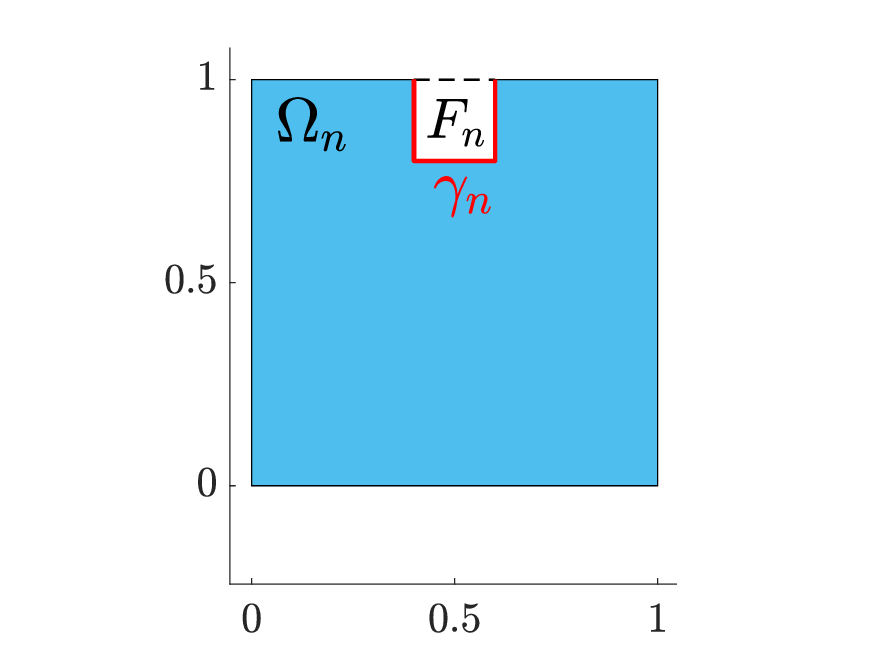}
			\caption{Geometry with negative feature.}
			\label{fig:negposfeat-a}
		\end{subfigure}\hspace{1cm}%
		\begin{subfigure}[t]{0.4\textwidth}
			\centering
			\includegraphics[width=0.96\linewidth]{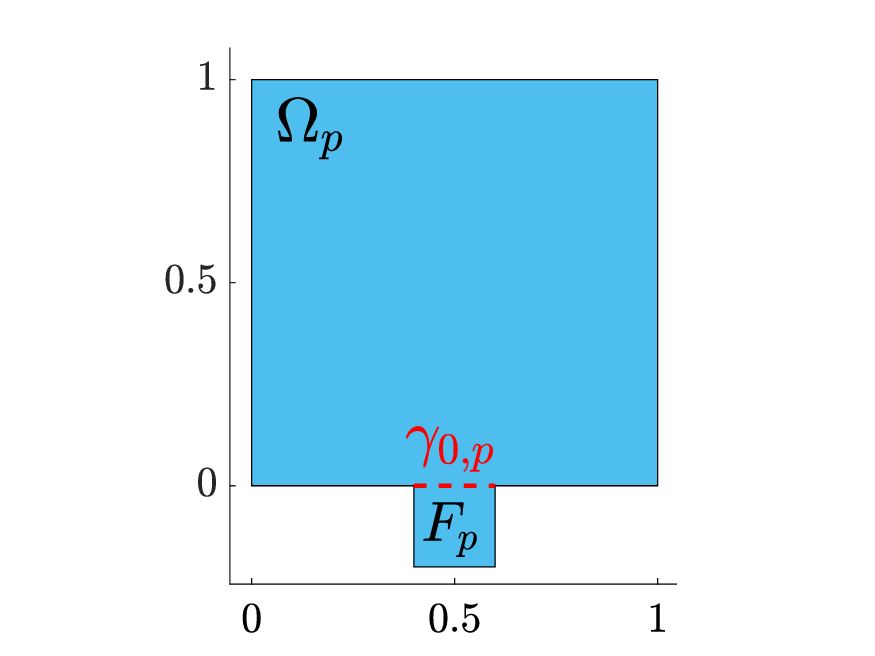}
			\caption{Geometry with positive feature.}
			\label{fig:negposfeat-b}
		\end{subfigure}	
		\caption{Test 2: Exact geometry with a negative and a positive feature.}
		\label{fig:negposfeat}
	\end{figure}
	
	In this second numerical example we consider the case of a positive and a negative boundary feature. As in the previous test case we choose $\Omega_0=(0,1)^2$, while we define
	$$F_\mathrm{n}=\left(\frac{1-\epsilon}{2},~\frac{1+\epsilon}{2}\right)\times \left(1-\epsilon,~1\right), $$ $$F_\mathrm{p}=\left(\frac{1-\epsilon}{2},~\frac{1+\epsilon}{2}\right)\times (-\epsilon,~0).$$
	For the negative feature case we chose as exact geometry $\Omega_n=\Omega_0\setminus \overline{F_\mathrm{n}}$, while for the positive feature case we choose $\Omega_\mathrm{n}=\text{int}(\overline{\Omega_0} \cup \overline{F_\mathrm{p}})$, as reported in Figure~\ref{fig:negposfeat}. 
	
	\noindent In both cases we consider the problem
	\begin{equation}
	\begin{cases}
	-\Delta u=f &\text{in } \Omega_\star\\
	u=0 &\text{on } \Gamma_{\mathrm{D}}\\
	\nabla u\cdot \bm{n}=0 &\text{on } \Gamma_{\mathrm{N}}
	\end{cases}\label{prob3}
	\end{equation}
	with $\star=\lbrace n,p \rbrace$, $f=1$, $\Gamma_{\mathrm{D}}=\lbrace (x,y): x=0 \vee x=1\rbrace$ and $\Gamma_{\mathrm{N}} =\partial \Omega_\epsilon\setminus \Gamma_{\mathrm{D}}$. We recall that, according to its definition, $\Gamma_{\mathrm{N}}$ includes also the feature boundary.
	
	\begin{figure}
		\centering
		\begin{subfigure}{0.4\textwidth}
			\centering
			\includegraphics[width=0.92\linewidth]{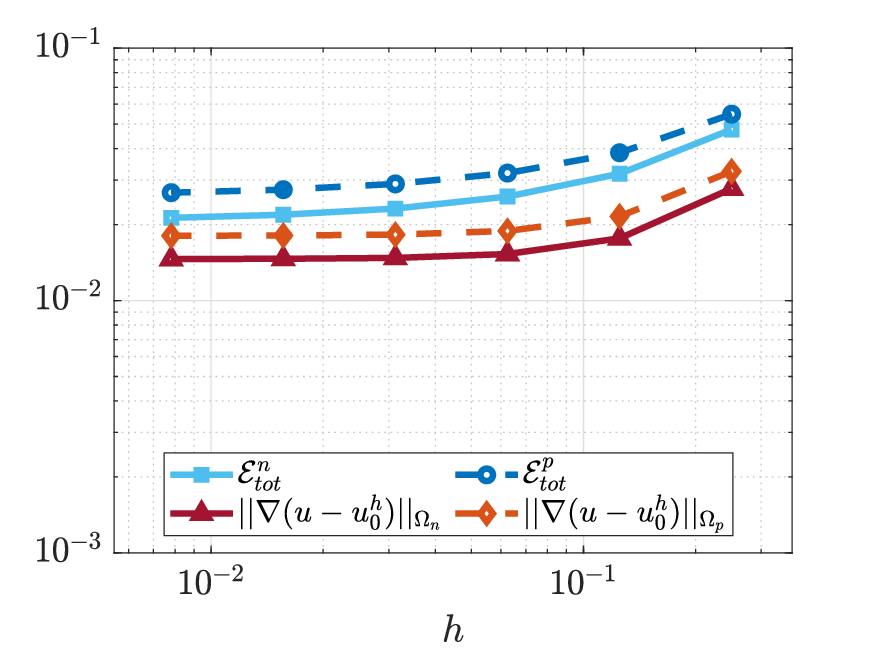}
			\caption{$\epsilon=0.2$.}
			\label{fig:test2_varH-a}
		\end{subfigure}\hspace{1cm}%
		\begin{subfigure}{0.4\textwidth}
			\centering
			\includegraphics[width=0.92\linewidth]{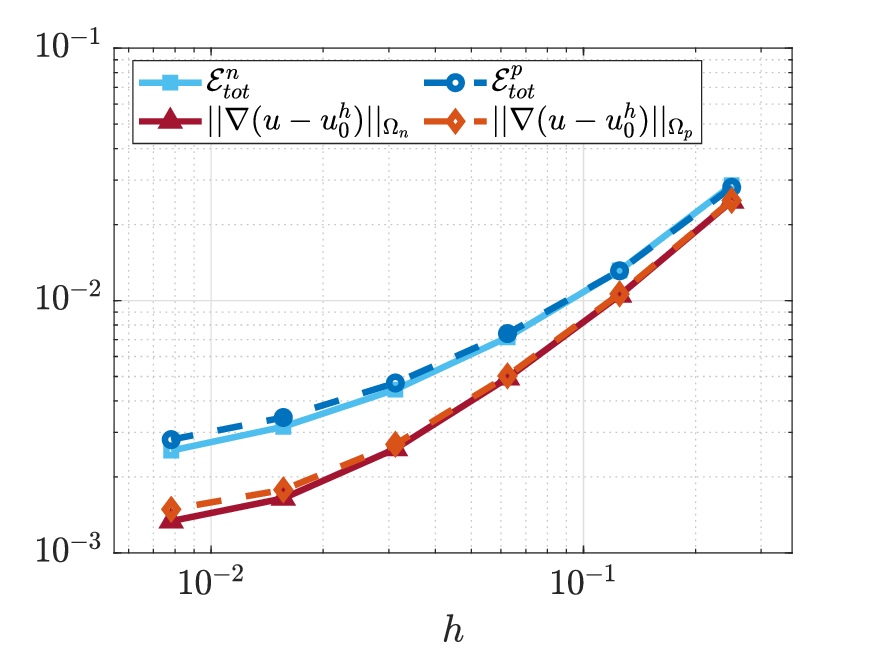}
			\caption{$\epsilon=0.05$.}
			\label{fig:test2_varH-b}
		\end{subfigure}
		\hfill
		\caption{Test 2: energy norm of the error and total estimator under mesh refinement, for the negative and the positive feature cases and for two different feature sizes.}
		\label{fig:test2_varH}
	\end{figure}
	\begin{figure}
		\centering
		\begin{subfigure}{0.4\textwidth}
			\centering
			\includegraphics[width=0.92\linewidth]{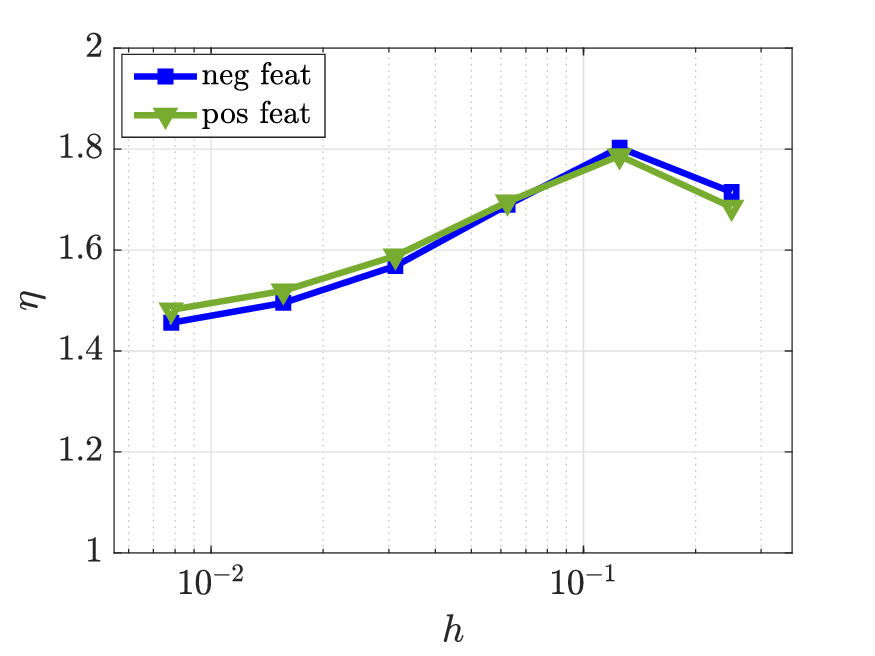}
			\caption{$\epsilon=0.2$}
			\label{fig:test2_effectivity-a}
		\end{subfigure}\hspace{1cm}%
		\begin{subfigure}{0.4\textwidth}
			\centering
			\includegraphics[width=0.92\linewidth]{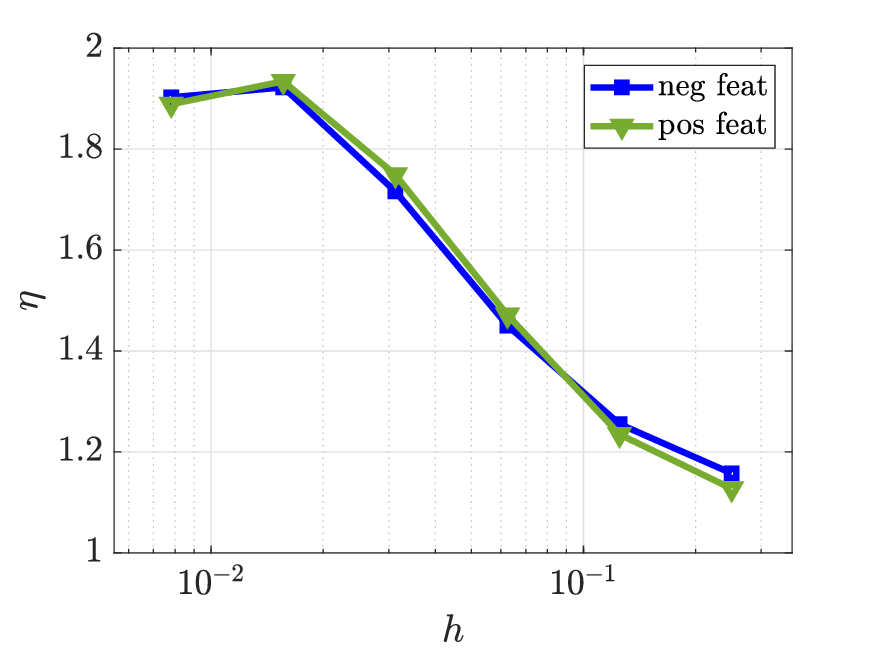}
			\caption{$\epsilon=0.05$}
			\label{fig:test2_effectivity-b}
		\end{subfigure}
		\caption{Test 2: effectivity index for the negative and the positive feature cases for two different feature sizes.}
		\label{fig:test2_effectivity}
	\end{figure}
	\begin{figure}
		\centering
		\begin{subfigure}[t]{.4\textwidth}
			\centering
			\includegraphics[width=0.92\linewidth]{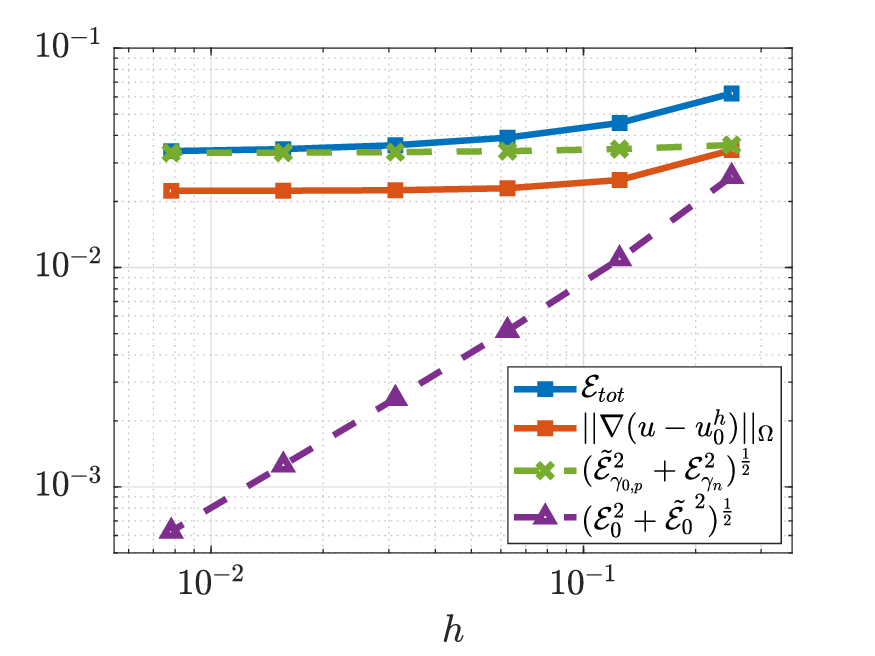}
			\caption{Energy norm of the error and total estimator with its components under mesh refinement. $\epsilon=0.2$}
			\label{fig:test2_2feat-a}
		\end{subfigure}\hspace{1cm}%	
		\begin{subfigure}[t]{.4\textwidth}
			\centering
			\includegraphics[width=0.92\linewidth]{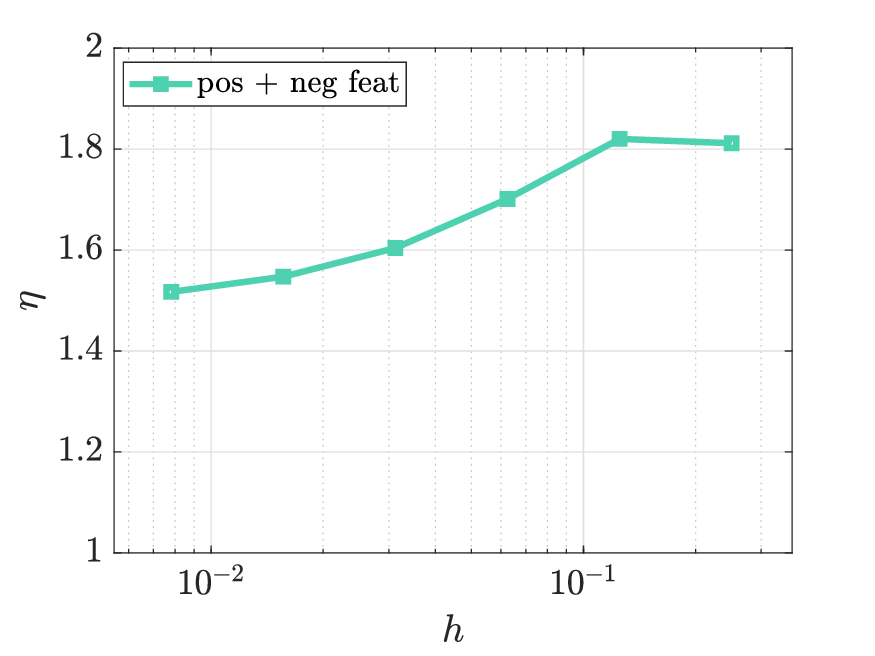}
			\caption{Effectivity index under mesh refinement. $\epsilon=0.2$}
			\label{fig:test2_2feat-b}
		\end{subfigure}
		\caption{Test 2, simultaneous presence of positive and negative feature:  energy norm of the error, total estimator with its components and effectivity index. For both features, $\epsilon=0.2$.}
		\label{fig:test2_2feat}
	\end{figure}
	
	Figure~\ref{fig:test2_varH} compares the convergence under mesh refinement of the total estimator and error in the negative and in the positive feature case and for two different values of $\epsilon$. Let us recall how, in presence of a negative feature the total estimator is defined as in \eqref{totestim_neg}, with $\estim{\gamma}=\estim{\gamma_\mathrm{n}}$ computed on $\gamma_\mathrm{n}$ (see Figure~\ref{fig:negposfeat-a}) from the equilibrated flux reconstructed on $\Omega_0$; in the case of a positive feature, instead, the definition is provided by \eqref{totestim_pos_F} with $\estimf{\gamma_0}=\estimf{\gamma_{0,\mathrm{p}}}$ computed on $\gamma_{0,\mathrm{p}}$ (see Figure~\ref{fig:negposfeat-b}) from the equilibrated flux reconstructed on the feature itself. In the following we choose $C_D=\tilde{C}_D=1$ and, for the sake of clarity, we denote respectively by $\estim{\mathrm{tot}}^\mathrm{n}$ and $\estim{\mathrm{tot}}^\mathrm{p}$ the total estimator in the negative and in the positive feature case. As in the previous test case, fixing the feature size and refining the mesh, we can observe how the overall error reaches a plateau, and how this behavior is captured also by the total estimator, both in the negative and in the positive feature case. As expected, a bigger feature (Figure~\ref{fig:test2_varH-a}) produces a stagnation of the error and of the estimator already for coarse meshes, while if the feature is smaller (Figure~\ref{fig:test2_varH-b}) the defeaturing source of error becomes relevant only for finer meshes. 
	
	Figure \ref{fig:test2_effectivity} shows the trend of the effectivity index related to the curves reported in Figure~\ref{fig:test2_varH}.  
	Both for the negative and the positive feature case we observe that, as expected, $\eta \sim 1$ when the numerical component is dominating (coarse meshes in Figure~\ref{fig:test2_varH-b}). We can instead observe how $\eta \sim 1.5$ when the defeaturing component dominates (fine meshes in Figure~\ref{fig:test2_varH-a}), and how the effectivity index is in general lower with respect to the internal negative feature case (Test 1), with $\eta<2$ even when both the defeaturing and the numerical component have a significant impact.

	Finally, Figure~\ref{fig:test2_2feat} refers to the case in which the negative and the positive features are simultaneously present, i.e. $\Omega=\text{int}(\overline{\Omega_0}\cup \overline{F_\mathrm{p}})\setminus \overline{F}_\mathrm{n}$. For both features we choose $\epsilon=0.2$. Figure \ref{fig:test2_2feat-a} reports the convergence of the error and of the estimator under mesh refinement. The total estimator is, in this case, defined as \begin{equation*}\estim{\mathrm{tot}}=C(\estimf{\gamma_{0,\mathrm{p}}}^2+\estim{\gamma_\mathrm{n}}^2)^\frac{1}{2}+(\estim{0}^2+\estimf{0}^2)^\frac{1}{2},
	\end{equation*}
	with $C>0$ being a constant independent of the size of both features. In particular, we choose here $C=1$. As in the previous test cases, the estimator appears to correctly capture the behavior of the overall error. The corresponding effectivity index is reported in Figure~\ref{fig:test2_2feat-b}.

	\subsection{Test 3: multiple internal features}
	
	For this last numerical example we consider a case with multiple internal features, similar to the one proposed in \cite{AC_2023,BCV2022_arxiv}. Our aim is to show the capability of the proposed estimator to identify the most relevant features and to provide a criterion to decide whether a feature should be added or not, according also to the magnitude of the numerical source of the error.
	
\begin{figure}
		\centering
		\begin{subfigure}[t]{0.4\textwidth}
			\centering
			\includegraphics[width=0.95\linewidth]{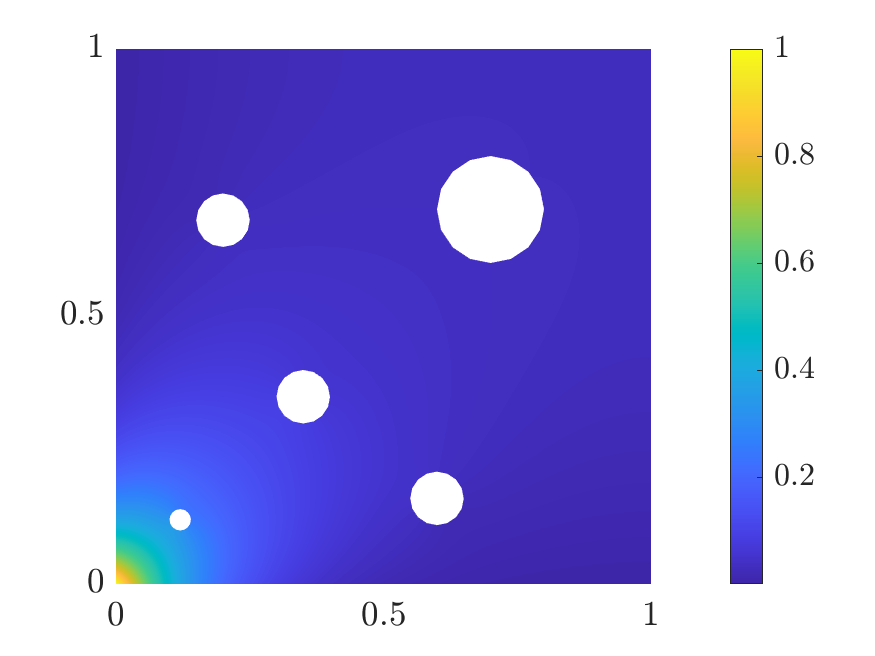}
			\caption{Reference solution.}
			\label{fig:sol_ex_test3-a}
		\end{subfigure}\hspace{1cm}%
		\begin{subfigure}[t]{0.4\textwidth}
			\centering
			\includegraphics[width=0.95\linewidth]{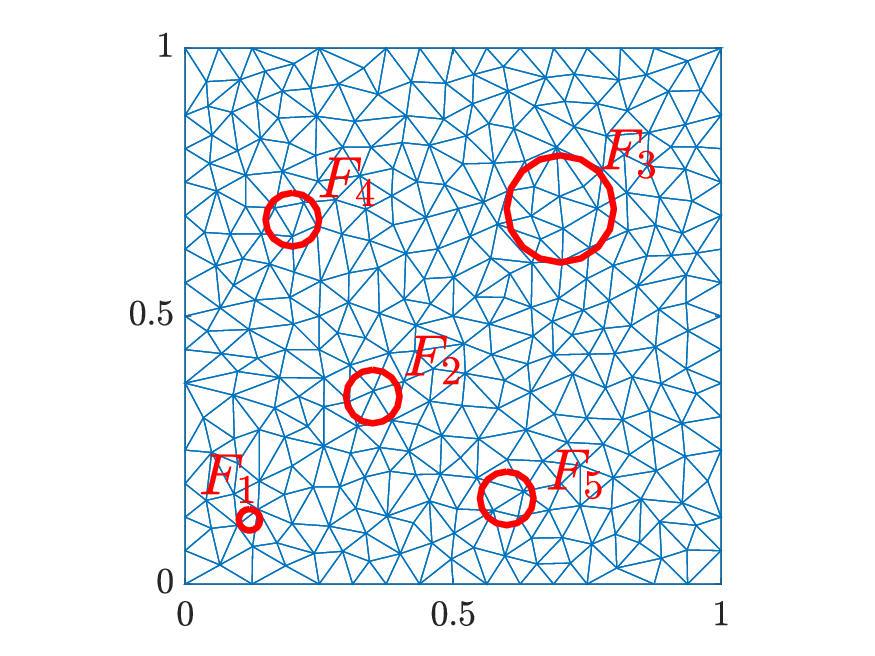}
			\caption{Example of computational mesh defined on the totally defeatured geometry $\Omega_0$.}
			\label{fig:sol_ex_test3-b}
		\end{subfigure}
		\caption{Test 3: reference solution $u$ and totally defeatured geometry $\Omega_0$ with an example of computational mesh.}
		\label{fig:sol_ex_test3}
	\end{figure}

	Let us define the defeatured geometry again as $\Omega_0=(0,1)^2$ and let us consider a set of $I$ features $\mathcal{F}=\lbrace{F_i}\rbrace_{i\in \mathcal{I}}$, $\mathcal{I}=\lbrace1,...,I\rbrace$, each of which is a polygon of 16 faces, inscribed in a circle of radius $\epsilon_i$ and centered in $\bm{x}_C^i$. In particular we choose $I=5$ and
	\begin{align*}
	&\bm{x}_C^1=(0.12,0.12), ~ \epsilon_1=0.02; \qquad 
	\bm{x}_C^2=(0.35,0.35),  ~ \epsilon_2=0.05\\
	&\bm{x}_C^3=(0.65,0.65), ~ \epsilon_3=0.10;\qquad 
	\bm{x}_C^4=(0.20,0.68),  ~ \epsilon_4=0.05;\\&
	\bm{x}_C^5=(0.65,0.16),  ~ \epsilon_5=0.05.
	\end{align*}
	The boundary of the $i$-th feature is denoted by $\gamma_i$. 
	On the exact geometry $\Omega=\Omega_0\setminus\overline{\bigcup_{i \in \mathcal{I}}F_i}$, which is reported in Figure \ref{fig:sol_ex_test3-a}, we consider the problem
	\begin{equation}
	\begin{cases}
	\Delta u=0 &\text{in } \Omega\\
	u=g_{\mathrm{D}} &\text{on } \Gamma_{\mathrm{D}}\\
	\nabla u\cdot \bm{n}=0 &\text{on } \Gamma_{\mathrm{N}}
	\end{cases}
	\end{equation}
	with $\Gamma_{\mathrm{D}}=\left\lbrace (x,y): x=0\vee y=0\right\rbrace, $
	$g_{\mathrm{D}}(x,y)=e^{-8(x+y)}$ and $\Gamma_{\mathrm{N}}=\partial \Omega \setminus \Gamma_{\mathrm{D}}$, including also the feature boundaries. The reference solution $u$ is also reported in Figure~\ref{fig:sol_ex_test3-a}.
	
	Let $\mathcal{M}\subseteq \mathcal{F}$ be a subset of features indexed by $j \in \mathcal{I}^\star\subseteq\mathcal{I}$ and let us denote by $\Omega_0^\mathcal{M}$ a generic partially defeatured geometry, obtained by including the features in $\mathcal{M}$ to the defeatured geometry, i.e.  $\Omega_0^\mathcal{M}=\Omega_0\setminus \overline{\bigcup_{j \in \mathcal{I}^\star}F_j}$. If $\mathcal{M}=\emptyset$, then $\Omega_0^\mathcal{M}=\Omega_0$. An example of a computational mesh defined on $\Omega_0$ is reported in Figure~\ref{fig:sol_ex_test3-b}. We will use $u_0^h$ to refer both to the numeric solution computed on $\Omega_0$ and to the one computed on $\Omega_0^\mathcal{M}$, the meaning being clear from the context.
	In presence of multiple negative features the total estimator is defined as $\estim{\mathrm{tot}}=\alpha_D\estim{\gamma}+\estim{0}$, where
	$$\estim{\gamma}=\Big(\sum_{k \in \mathcal{I} \setminus \mathcal{I}^\star }\estim{\gamma_k}^2\Big)^\frac{1}{2}.$$ 
	
	\begin{table}
		\renewcommand*{\arraystretch}{1.45}
		\centering 
		\caption{Test 3: components of the total estimator for differently refined mesh and for different choices of the (partially) defeatured geometry.}
		\label{test3_table}
		\begin{tabular}{c|c|cc|ccccc|c}
			\hline
			$\bm{h}$&$\bm{\mathcal{M}}$&$\bm{N_{dof}}$&	$\estim{0}$ &$\estim{\gamma_1}$ &$\estim{\gamma_2}$ &$\estim{\gamma_3}$ &$\estim{\gamma_4}$ &$\estim{\gamma_5}$&$\estim{\gamma}$ \\
			\hline
			\multirow{2}{*}{$6.25e-2$}&$\bm{\emptyset}$&  1240&	$0.101$&${0.147}$ &$0.050$ &$0.008$&$0.026$ &$0.036$&0.162\\
			&$\bm{F_1}$&  1380&	$0.096$&$-$ &$0.048$ &$0.008$&$0.025$ &$0.035$&0.065\\
			\hline\hline
			\multirow{3}{*}{$3.13e-2$}&$\bm{\emptyset}$&  4960&	$0.052$&${0.146}$ &$0.050$ &$0.008$&$0.025$ &$0.036$&0.162\\
			&$\bm{F_1}$&  5501&	$0.050$&$-$ &$0.048$ &$0.008$&$0.025$ &$0.035$&0.065\\
			&$\bm{F_1,F_2}$&  5684&	$0.053$&$-$ &$-$ &$0.007$&$0.024$ &$0.035$&0.043\\
			\hline\hline
			\multirow{5}{*}{$1.56e-2$}&$\bm{\emptyset}$&  19840&	$0.027$&${0.146}$ &$0.050$ &$0.008$&$0.025$ &$0.036$&0.161\\
			&$\bm{F_1}$&  21963&	$0.026$&$0$ &$0.048$ &$0.008$&$0.025$ &$0.035$&0.065\\
			&$\bm{F_1,F_2}$&  22618&	$0.027$&$-$ &$-$ &$0.007$&$0.024$ &$0.034$&0.042\\
			&$\bm{F_1,F_2,F_5}$&  23953&	$0.027$&$-$ &$-$ &$0.007$&$0.024$ &$-$&0.025\\
			&$\bm{F_1,F_2,F_4,F_5}$&  24372&	$0.027$&$-$ &$-$ &$0.007$&$-$ &$-$&0.007
		\end{tabular}
	\end{table}
	
	Table \ref{test3_table} reports the value of the components of the total estimator for differently refined meshes and for different choices of $\mathcal{M}$, i.e. of the partially defeatured geometry on which $u_0^h$ is computed. The rows of the table are divided into three sets, corresponding to three differently refined meshes. The variation in the number of degrees of freedom which can be observed when a feature is included into the geometry is related to the adaptation of the mesh to the feature boundary and to the deletion of the degrees of freedom lying inside the feature itself. Looking at the columns from 5 to 9 we can see how feature $F_1$ is clearly the most relevant, since $\estim{\gamma_1}>\estim{\gamma_i}$ for all $i>1$. 
	This is expected since, although being the smallest feature, it is located in a region in which the gradient of the solution is very steep. Feature $F_3$ is, instead, almost irrelevant: despite being the biggest one it is located in a region in which the solution is rather flat and hence its impact on the solution accuracy tends to be negligible. As expected, the values of $\estim{\gamma_i}$ are independent from the mesh size, meaning that the relevance of the features can be evaluated even on a coarse mesh. However, the choice of including the $i$-th feature into the geometry should be taken by comparing $\estim{\gamma_i}$ with $\estim{0}$, which is a sharp indicator of the numerical source of error. In particular, a value of $\estim{\gamma_i}$ considerably bigger than $\estim{0}$ means that we will not be able to significantly reduce the error by mesh refinement, unless the feature is added. This is true, for example, for feature $F_1$ with the second considered mesh and for features $F_1$ and $F_2$ for the finest mesh.
	
	Table \ref{test3_table2} focuses exactly on these cases, reporting the values of the energy norm of the overall error and of the total estimator, along with the corresponding effectivity index. In particular, $\estim{\mathrm{tot}}$ is computed with $\alpha_D=1$. In Table \ref{test3_table}, we can observe that, if $\mathcal{M}=\emptyset$, the reduction of $\estim{0}$ when going from $h=3.13\cdot 10^{-2}$ ($\sim5k$ degrees of freedom) to $h=1.56\cdot 10^{-2}$ ($\sim20k$ degrees of freedom) is of about $50\%$. This is expected, since $\estim{0}$ should converge as $\mathcal{O}(h)$. However, looking at rows 1 and 3 in Table \ref{test3_table2}, we see that the drop in the total estimator, and hence in the error, is under $20\%$. Adding feature $F_1$ and refining the mesh at the same time the drop is instead of about $60\%$, as it can be seen by comparing rows 1 and 4 in Table \ref{test3_table2}. Let us remark that, for the finest considered mesh, also feature $F_2$ becomes rather relevant. However, $\estim{\gamma_2}$ is closer to $\estim{0}$ and hence adding it to the geometry has a smaller impact on the solution accuracy.
	\begin{table}
		\renewcommand*{\arraystretch}{1.45}
		\centering 
		\caption{Test 3: energy norm of the error, total estimator and effectivity index for differently refined mesh and for different choices of the (partially) defeatured geometry. }
		\label{test3_table2}
		\begin{tabular}{c|c|c|cc}
			\hline
			$\bm{h}$&$\bm{\mathcal{M}}$ &$||\nabla(u-u_0^h)||_\Omega$&$\estim{\mathrm{tot}}$ & $\eta~$
			\\
			\hline
			\multirow{2}{*}{$3.13e-2$}&$\bm{\emptyset}$&$0.079$&$0.214$&2.71
			\\
			&$\bm{F_1}$&$0.052$&$0.115$&2.21
			\\
			\hline\hline
			\multirow{3}{*}{$1.56e-2$}&$\bm{\emptyset}$&$0.066$&$0.188$&2.85
			\\
			&$\bm{F_1}$&$0.033$&$0.091 $&2.76
			\\
			&$\bm{F_1,F_2}$&$0.028$&$0.069$&2.46
		\end{tabular}
	\end{table}
	
	This experiment is to be intended as a preliminary test for the use of the proposed estimator in an adaptive strategy, involving both geometrical adaptation (i.e. feature inclusion) and local mesh refinement. We decide to leave this to a forthcoming work: indeed the procedure adopted for the computation of the equilibrated flux requires the mesh to be conforming to the domain boundaries, and hence also to the boundary of the features which are actually included in the partially defeatured geometry. However, to build an efficient and flexible adaptive algorithm we do not want to remesh the geometry each time a feature is added, and for this reason a generalization of the equilibrated flux reconstruction to trimmed meshes needs to be considered.

	\section*{Conclusions}
	In this work we have proposed a new \emph{a posteriori} error estimator for defeaturing problems based on an equilibrated flux reconstruction and designed for finite elements. The Poisson equation with Neumann boundary conditions on the feature boundary was taken as a model problem. The reliability of the estimator has been proven both in the negative and in the positive feature case, and tested with several numerical examples. The choice of using an equilibrated flux reconstruction leads to an estimator which is able to bound sharply the numerical component of the error and which never requires to evaluate the normal trace of the numerical flux, which is typically discontinuous on element edges in a standard finite element discretization.
	
	This work is to be intended as a preliminary analysis for the use of the proposed estimator in an adaptive strategy, allowing not only for mesh refinement, but also for an automatic inclusion of those features whose absence causes most of the accuracy loss. The proposed estimator does not require the mesh to be conforming to the feature boundary until the feature is included in the computational domain itself. Indeed, computing the integral of the normal trace of the equilibrated flux reconstruction on a generic curve is always possible, regardless of the intersections with the mesh elements. However, the procedure which was adopted to reconstruct the equilibrated flux is designed for meshes which are conforming to the computational domain boundary and this would require to remesh the domain each time a feature is added by the adaptive procedure, hence increasing the complexity of the algorithm. For this reason, an extension of the equilibrated flux reconstruction to the case of trimmed meshes needs to be considered, so that the geometry never needs to be remeshed. This generalization is left to a forthcoming work, that is currently under preparation. Although the proof of the reliability of the estimator holds in $\mathbb{R}^d$, $d=2,3$, and for any polynomial order $p\geq 1$, we decided to propose numerical experiments only in $\mathbb{R}^2$ and for $p=1$. The application of the estimator on more complex, realistic and tridimensional geometries, or the use of a higher order finite element approximation, are left to a forthcoming work as well, both extensions having an impact only on implementation aspects.

\section*{Acknowledgments}
The authors are grateful to Zhaonan Dong (Inria, Paris) for sharing his code on equilibrated fluxes and to Paolo Bardella (DET, Politecnico di Torino) for the MeshToolbox library.\smallskip

The authors acnowledge the support of the Swiss National Science Foundation (via project MINT n. 200021\_215099, PDE tools for analysis-aware geometry processing in simulation science) and of
European Union Horizon 2020 FET program (under grant agreement No 862025 (ADAM2)).

Author O. Chanon acknowledges the support of the Swiss National Science Foundation through the project n.P500PT 210974.

\printbibliography
\end{document}